\documentstyle[12pt,times,graphicx,subfigure,amssymb]{article}
\newcommand\beq[1]{ \begin{equation}\label{#1} }
\newcommand{\eeq}{ \end{equation} }
\newcommand{\beqno}{ \[ }
\newcommand{\eeqno}{ \] }
\newcommand\beqa[1]{ \begin{eqnarray} \label{#1}}
\newcommand{\eeqa}{ \end{eqnarray} }
\newcommand{\beqano}{ \begin{eqnarray*} }
\newcommand{\eeqano}{ \end{eqnarray*} }
\newcommand\equ[1]{{\rm (\ref{#1})}}
\newcommand{\eps}{\varepsilon}
\newcommand{\R}{\mathbb{R}}
\newcommand{\Z}{\mathbb{Z}}

\newcommand{\T}{\mathbb{T}}
\newcommand{\C}{\mathbb{C}}

\begin{document}

\title{Stability of nearly--integrable systems with dissipation}

\author{ Alessandra Celletti\\
{\footnotesize Dipartimento di Matematica}\\
{\footnotesize Universit\`a di Roma Tor Vergata}\\
{\footnotesize Via della Ricerca Scientifica 1}\\
{\footnotesize I-00133 Roma (Italy)}\\
{\footnotesize \texttt{(celletti@mat.uniroma2.it)}}
\and
Christoph Lhotka\\
{\footnotesize Dipartimento di Matematica}\\
{\footnotesize Universit\`a di Roma Tor Vergata}\\
{\footnotesize Via della Ricerca Scientifica 1}\\
{\footnotesize I-00133 Roma (Italy)}\\
{\footnotesize \texttt{(lhotka@mat.uniroma2.it)}}
}

\maketitle

\vglue1cm

\begin{abstract}
We study the stability of a vector field associated to a
nearly--integrable Hamiltonian dynamical system to which a
dissipation is added. Such a system is governed by two
parameters, named the perturbing and dissipative parameters, and
it depends on a drift function. Assuming that the frequency of
motion satisfies some resonance assumption, we investigate the
stability of the dynamics, and precisely the variation of the
action variables associated to the conservative model. According
to the structure of the vector field, one can find linear and
exponential stability times, which are established under smallness
conditions on the parameters. We also provide some applications to
concrete examples, which exhibit a linear or exponential
stability behavior.
\end{abstract}

\vglue1cm


\vglue1cm

\tableofcontents

\section{Introduction}
We investigate the behavior of nearly--integrable Hamiltonian
vector fields to which a dissipative contribution is added. The
vector field is ruled by two parameters, namely the perturbing
parameter (measuring the non--integrability of the system) and the
dissipative parameter (providing the size of the dissipative
term). We assume that the phase space is contracted by time evolution.
A drift
function enters the equations of motion as an unknown function;
its role is fundamental, since it must be properly chosen in order
to meet some compatibility conditions ensuring the existence of a
normal form (compare with KAM results as in \cite{CCL}). We concentrate on the behavior of the variables which
are actions of the conservative system (i.e. setting to zero the
dissipative parameter). We assume that the initial conditions
define a resonant frequency for the integrable conservative system
(i.e. setting to zero both the perturbing and dissipative
parameters). Under smallness conditions on the parameters, we
prove that the action variables stay locally bounded over a given
time interval (see also \cite{nekh}). The length of the time interval depends on the
functions defining the equations of motion and, precisely, whether
there appear also dissipative resonant terms in the original as
well as in the normalized vector field. Notice that such a result
provides a useful information concerning the transient time,
namely the time needed to reach the attractor.

\noindent The proof of the result is based on the construction of
a suitable coordinate transformation, which is provided by the
composition of a conservative and a dissipative change of
variables. A similar technique, based on a non--resonant normal
form, has been already implemented in \cite{cellho2010a} in order
to investigate a vector field of the type studied in this paper, but
in the simplest case of non--resonant frequency. In this case, under smallness conditions on the
parameters, one can prove that the actions stay always bounded for
exponential times. As in classical perturbation theory, the first
transformation removes the conservative perturbation to higher
orders (see, e.g., \cite{CE10}); the corresponding normal form is
composed by resonant or average terms. The second transformation
is performed to normalize the dissipative terms; the normal form
equations defining the dissipative change of variables can be
solved, provided that the drift function is chosen in such a way
that the compatibility condition is satisfied. The final normal
form contains just resonant and average terms up to a given order
in the perturbing and dissipative parameters. As in classical
Nekhoroshev's theorem (\cite{nekh}, \cite{nekh1}, see also
\cite{benetal}, \cite{gionekh}) by properly choosing the order of
the normal form, one can determine stability bounds. The stability
time is exponential, whenever conservative resonant terms do not
appear in the equation for the time variation of the normalized
action variable or whenever the dissipative resonant contributions
are zero. In the other cases the stability time depends on the
inverse of the product of the perturbing and dissipative
parameters. The scheme of the proof, which is presented for a
non--autonomous, time--periodic system (see also \cite{giozeh}) follows
closely \cite{poe1}, where a very clear and enlightening proof of
Nekhoroshev's theorem is given. The proof is constructive and it
allows us to provide explicit expressions for the conservative and
dissipative transformations (see also \cite{BM}, \cite{DGL},
\cite{LI05}). In our opinion there are several physical problems,
which can be analyzed by our method. For example, there are many
results concerning the stability of the (resonant) Lagrangian
points in a conservative framework (see \cite{CF}, \cite{CG},
\cite{GDFGS}, \cite{GS}, \cite{LED}), but none of them takes
into account dissipative effects (like Solar radiation,
Poynting--Roberston drag, Yarkowsky effect, etc.), which might
significantly affect the dynamics. In this respect, we believe that it would
be interesting to analyze these models including a dissipative effect
by using the results contained in this paper. We provide examples
of normal forms in some concrete one--dimensional, time--dependent
model problems, which illustrate different cases corresponding to
linear (i.e., proportional to the inverse of the product of the
perturbing and dissipative parameters) or exponential stability
times. We also provide an application of the theorem in order to
obtain rigorous stability bounds for the previous model problems.

\vskip.1in

\noindent
The paper is organized as follows. Notations and assumptions are defined in Section~\ref{sec:notation}.
The resonant normal form Lemma and the stability Theorem are proven in Section~\ref{sec:stability}.
Examples of normal form constructions to concrete model problems is given in Section~\ref{sec:examples}.
An application of the stability theorem is provided in Section~\ref{sec:APP}.

\section{Notations and assumptions}\label{sec:notation}
We introduce the $\ell$--dimensional, time--dependent vector field, described by the equations
\beqa{SE}
\dot x&=&\omega(y)+\varepsilon h_{10,y}(y,x,t)+\mu f_{01}(y,x,t)\nonumber\\
\dot y&=&-\varepsilon h_{10,x}(y,x,t)-\mu \Big(g_{01}(y,x,t)-
\eta(y,x,t)\Big)\ , \eeqa where $y\in\R^\ell$,
$(x,t)\in{\T}^{\ell+1}$, while the definitions and assumptions on
the parameters and functions are the following\footnote{The
subscripts $x$, $y$, $t$ denote derivatives with respect to $x$,
$y$, $t$, i.e. $h_x\equiv {{\partial h}\over {\partial x}}$,
$h_y\equiv {{\partial h}\over {\partial y}}$, $h_t\equiv
{{\partial h}\over {\partial t}}$.}.

\begin{enumerate}
    \item Having fixed an initial datum
$y_0\in{\R}^\ell$, we denote by $A\subset {\R}^\ell$ an open
neighborhood of $y_0$.
    \item The vector field depends on the parameters $\varepsilon\in{\R}_+$ (perturbing parameter), $\mu\in{\R}_+$ (dissipative parameter);
    we remark that we could equally admit vector parameters, i.e. $\varepsilon\in{\R}_+^\ell$, $\mu\in{\R}_+^\ell$, but for simplicity of exposition we present the details just for the scalar case $\varepsilon\in{\R}_+$, $\mu\in{\R}_+$.
    \item The functions $\omega$ and $\eta$ are
real--analytic, $\ell$--dimensional vector functions with
components $(\omega^{(1)},...,\omega^{(\ell)})$ and
$(\eta^{(1)},...,\eta^{(\ell)})$. We assume that there exists a
regular function $h_{00}=h_{00}(y)$ such that ${{\partial
h_{00}(y)}\over {\partial y}}=\omega(y)$. Let
$h_0(y,u)=h_{00}(y)+u$, $u\in{\R}$, be the unperturbed Hamiltonian
function associated to the conservative vector field ($\mu=0$) in
the extended phase space. Let $\omega_e(y)\equiv(\omega(y),1)$ be
the frequency vector in the extended phase space. Following
\cite{poe1} we make the hypothesis that $h_{0}$ is $L,M$--quasi
convex, namely there exist $L$, $M>0$, such that for all
$z\equiv(y,u)\in A\times{\R}$ at least one of the following
inequalities is satisfied: \beq{QC} |\omega_e(y)\cdot v|>L\|v\|\
,\qquad {{\partial^2 h_{0}(z)}\over {\partial z^2}}v\cdot v\geq
M\|v\|^2\ ,\qquad \forall v\in{\R}^{\ell+1}\ , \eeq where the dot
denotes the scalar product and $\|\cdot\|$ denotes the Euclidean
norm.
    \item In the following we will use the
    vector field \equ{SE} in the extended phase space with $\dot t=1$ and with $u$ conjugated to time:
\beqa{SEexp}
\dot x&=&\omega(y)+\varepsilon h_{10,y}(y,x,t)+\mu f_{01}(y,x,t)\nonumber\\
\dot y&=&-\varepsilon h_{10,x}(y,x,t)-\mu \Big(g_{01}(y,x,t)-
\eta(y,x,t)\Big)\nonumber\\
\dot u&=&-\varepsilon h_{10,t}(y,x,t)+\mu \sigma(y,x,t)\ , \eeqa
where the unknown function $\sigma$ is introduced for later
convenience (see next point).
    \item We assume that $f_{01}$, $g_{01}$, $\eta$
are real--analytic, $\ell$--dimensional
vector functions from $A\times {\T}^{\ell+1}$ to ${\R}^\ell$,
while $h_{10}$, $\sigma$ are periodic and real--analytic from
$A\times {\T}^{\ell+1}$ to ${\R}$. We remark that $\eta$, $\sigma$
are unknown functions, which will be properly chosen so to meet
some compatibility requirements in order to obtain a suitable
normal form (see Section~\ref{sec:stability}).
    \item We assume that the vector
field is dissipative and that the phase space volume is contracted
by the time evolution.
    \item For a given initial datum $y_0=y(0)\in A$, we assume that there exists a lattice $\Lambda\subset{\Z}^{\ell+1}$, such that
     the vector function $\omega=\omega(y_0)$ satisfies the \sl resonance condition \rm
\beq{res}
|\omega(y_0)\cdot k+j|=0\quad {\rm for\
all\ } (k,j)\in\Lambda\ .
\eeq
We also assume that there exists $K\in{\Z_+}$, $a>0$ and a subset $D\subseteq A$, such that for any $y\in D$
the following condition is satisfied:
\beq{omega0}
|\omega(y)\cdot k+j|\geq a\qquad {\rm for\ all}\ (k,j)\in{\Z}^{\ell+1}\backslash\Lambda\ ,\
|k|+|j|\leq K\ ,
\eeq
where for $k=(k_1,...,k_\ell)\in{\Z}^\ell$ we define the norm $|k|\equiv |k_1|+...+|k_\ell|$.
\item  We refer to $\eta=\eta(y,x,t)$ as the \sl drift \rm vector function with components
$(\eta^{(1)}(y,x,t)$, $...$, $\eta^{(\ell)}(y,x,t))$ that we expand as
$$
\eta^{(k)}(y,x,t)=\sum_{m=0}^\infty \sum_{j=0}^m
\eta_{j,m-j+1}^{(k)}(y,x,t) \varepsilon^j\mu^{m-j}\ ,\qquad
k=1,...,\ell\ .
$$
In a similar way we expand $\sigma$ in \equ{SEexp} as
$$
\sigma(y,x,t)=\sum_{m=0}^\infty \sum_{j=0}^m
\sigma_{j,m-j+1}(y,x,t) \varepsilon^j\mu^{m-j}\ .
$$

   \end{enumerate}

\noindent
\bf Remark 1. \rm We remark that for $\mu=0$ the equations \equ{SEexp} reduce to the conservative vector field, associated
to the nearly--integrable Hamiltonian function in the extended phase space
\beq{H}
{\cal H}(y,x,u,t)=h_{00}(y)+u+\varepsilon h_{10}(y,x,t)\ ,
\eeq
where $\omega(y)={\partial h}_{00}(y)/\partial y$.
Notice that the Hamiltonian \equ{H} is integrable
as far as the perturbing parameter is zero, i.e. $\varepsilon=0$. Since the vector
field \equ{SEexp} is dissipative (see assumption 6.), the energy associated to
\equ{H} is decreasing with time.

\vskip.1in

\noindent
We adopt the following notations and definitions for functions and norms.
\begin{description}
    \item $i)$ Integer subscripts denote the order in the
    perturbing and dissipative parameters, i.e.
    $F_{ij}$ denotes a function of order $\varepsilon^i\mu^j$.
\item $ii)$ For a function $f=f(y,x,t)$ and for any positive integer
$K$, we decompose $f$ as
$$
f(y,x,t)=\bar f(y)+f^{(nr,\leq K)}(y,x,t)+f^{(r,\leq K)}(y,x,t)+f^{(>K)}(y,x,t)\ ,
$$
being, respectively,
the average, the sum over the non--resonant
components with Fourier modes less or equal than $K$,
the projection over the resonant space defined by the lattice $\Lambda$ excluding the origin with Fourier modes less or equal than $K$,
the sum over the components with Fourier modes greater than $K$, namely:
\beqano
\bar f(y)&\equiv&
{1\over {(2\pi)^{\ell+1}}}\int_{{\T}^{\ell+1}} f(y,x,t)\ dx dt\nonumber\\
f^{(nr,\leq K)}(y,x,t)&\equiv& \sum_{(k,j)\in{\Z}^{\ell+1}\backslash\Lambda,\ |k|+|j|\leq K}
\hat f_{kj}(y)e^{\imath(k\cdot x+jt)}\nonumber\\
f^{(r,\leq K)}(y,x,t)&\equiv& \sum_{(k,j)\in\Lambda\backslash\{0\},\ |k|+|j|\leq K}
\hat f_{kj}(y)e^{\imath(k\cdot x+jt)}\nonumber\\
f^{(>K)}(y,x,t)&\equiv& \sum_{(k,j)\in{\Z}^{\ell+1},\ |k|+|j|> K}
\hat f_{kj}(y)e^{\imath(k\cdot x+jt)}\ , \eeqano where
$\imath=\sqrt{-1}$ and $\hat f_{kj}$ are the Fourier coefficients.
    \item $iii)$ We say that a function is of order $k$ in $\varepsilon$ and $\mu$, in symbols
$O_k(\varepsilon,\mu)$, if its Taylor series expansion
in $\varepsilon$, $\mu$ contains powers of $\varepsilon^i\mu^j$ with
$i+j\geq k$, $i\geq 0$, $j\geq 0$.
\item $iv)$ We denote by $C_{r_0}(A)$ the complex neighborhood of $A$
of radius $r_0$, namely
$$
C_{r_0}(A)\equiv \{{y\in\C^\ell:}\ \|y-y_A\|\leq r_0\ {\rm for\ all\ } y_A\in A\}\ .
$$
Moreover, let $C_{s_0}({\T}^{\ell+1})$ be
the complex strip of radius $s_0$ around ${\T}^{\ell+1}$, namely
$$
C_{s_0}({\T}^{\ell+1})\equiv \{(x,t)\in{\C}^{\ell+1}:\ \max_{1\leq j\leq \ell}|\Im(x_j)|
\leq s_0\ ,\ |\Im(t)|\leq s_0\}\ ,
$$
where $\Im$ denotes the imaginary part.
\item $v)$ Denoting the Fourier expansion of a function $f=f(y,x,t)$ as
$$
f(y,x,t)=\sum_{(k,j)\in{\Z}^{\ell+1}} \hat f_{kj}(y) e^{\imath(k\cdot x+jt)}\ ,
$$
we introduce the norm
$$
\|f\|_{r_0,s_0}\equiv \sup_{y\in C_{r_0}(A)}\ \sum_{(k,j)\in{\Z}^{\ell+1}}
|\hat f_{kj}(y)| e^{(|k|+|j|)s_0}\ .
$$
For a function $g=g(y)$ we define $\|g\|_{r_0}\equiv \sup_{y\in C_{r_0}(A)} \|g(y)\|$,
where $\|\cdot\|$ denotes the Euclidean norm. For a vector function $f=(f_1,..,f_\ell)$ we define
$$
\|f\|_{r_0,s_0}\equiv \sqrt{\sum_{j=1}^\ell \|f_j\|_{r_0,s_0}^2}\ .
$$
   \end{description}

\section{Bounds on the variation of the action variables}\label{sec:stability}
In order to bound the variation of the action variables, we
implement a change of coordinates such that the vector field
\equ{SEexp} is transformed to a resonant normal form, up to a suitable order $N$.
To this end we introduce a change of coordinates close to the identity and leaving time unaltered, say
\beq{c}
(Y,X,U,t)=\Xi^{(N)}(y,x,u,t)\ ,\qquad (Y,U)\in{\R}^{\ell+1}\ , \quad
(X,t)\in{\T}^{\ell+1}\ ,
\eeq
where $\Xi^{(N)}$ depends parametrically also on
$\varepsilon$, $\mu$, $\Xi^{(N)}=\Xi^{(N)}(y,x,u,t;\varepsilon,\mu)$ with
$\Xi^{(N)}(y,x,u,t;0,0)=Id$.
Let $K$ be as in \equ{omega0}; in the forthcoming Resonant Normal Form Lemma we aim to determine the
transformation of coordinates \equ{c}, so that \equ{SEexp} takes a resonant normal form of order $N$, that we write as
\beqa{NF}
\dot X&=&\Omega_d^{(N)}(Y)+F^{(r,\leq K)}(Y,X,t)+
F_{N+1}(Y,X,t)+F^{(>K)}(Y,X,t)\nonumber\\
\dot Y&=&G^{(r,\leq K)}(Y,X,t)+G_{N+1}(Y,X,t)+G^{(>K)}(Y,X,t)\nonumber\\
\dot U&=&H^{(r,\leq K)}(Y,X,t)+H_{N+1}(Y,X,t)+H^{(>K)}(Y,X,t)\ ,
\eeqa
where $\Omega_d^{(N)}:\R^\ell\rightarrow\R^\ell$ is the normalized frequency, related to
$\omega(Y)$ by
$$
\Omega_d^{(N)}(Y)=\Omega_d^{(N)}(Y;\varepsilon,\mu)\equiv
\omega(Y)+\sum_{i=1}^N\Omega_{i0}(Y)\varepsilon^i+
\sum_{j=1}^N\sum_{i=0}^{N-j}\Omega_{ij}(Y)\varepsilon^i\mu^j\ ,
$$
where $\Omega_{ij}(Y)$ are known vector functions;
$F^{(r,\leq K)}$, $G^{(r,\leq K)}$, $H^{(r,\leq K)}$ have Fourier components
belonging to the resonant lattice $\Lambda\backslash\{0\}$ with $F^{(r,\leq K)}$ depending
on both $\eps$, $\mu$, while $G^{(r,\leq K)}$, $H^{(r,\leq K)}$ depend only on $\eps$; $F_{N+1}$, $G_{N+1}$, $H_{N+1}$ are vector functions of order $O_{N+1}(\varepsilon,\mu)$;
$F^{(>K)}$, $G^{(>K)}$, $H^{(>K)}$ denote functions with Fourier modes greater than $K$.

\vskip.1in

\noindent
Similarly to \cite{cellho2010a} we decompose the coordinate transformation $\Xi^{(N)}$ as the composition
of two transformations $\Xi^{(N)}_c$ (conservative part) and $\Xi^{(N)}_d$ (dissipative part):
\beq{ast1}
(Y,X,U,t)=\Xi^{(N)}_d \circ \Xi^{(N)}_c(y,x,u,t)\ .
\eeq
Setting $(\tilde y,\tilde x,\tilde u,t)\equiv \Xi^{(N)}_c(y,x,u,t)$, the conservative
transformation $\Xi^{(N)}_c$ is defined through a
sequence of generating functions close to the identity, say $\psi_{j0}=\psi_{j0}(\tilde
y,x,t)$, $j=1,...,N$, such that
\beqa{A}
\tilde x&=&x+\sum_{j=1}^N \psi_{j0,y}(\tilde y,x,t)\varepsilon^j\equiv x+\psi_y^{(N)}(\tilde y,x,t)\nonumber\\
y&=&\tilde y+\sum_{j=1}^N \psi_{j0,x}(\tilde
y,x,t)\varepsilon^j\equiv \tilde y+\psi_x^{(N)}(\tilde y,x,t)\nonumber\\
u&=&\tilde u+\sum_{j=1}^N \psi_{j0,t}(\tilde
y,x,t)\varepsilon^j\equiv \tilde u+\psi_t^{(N)}(\tilde y,x,t)\ .
\eeqa
Notice that we can assume that the functions $\psi_{j0}$ (as well as $\alpha_{jk}$, $\beta_{jk}$, $\gamma_{jk}$  in \equ{B} below) do not depend on $u$, since
the functions appearing in \equ{SE} (or equivalently in \equ{SEexp}) do not depend on $u$.
We denote the inversion of \equ{A} as
\beqa{inv}
x&=&x(\tilde y,\tilde x,t)=\tilde x+\Gamma^{(x,N)}(\tilde y,\tilde x,t)\nonumber\\
y&=&y(\tilde y,\tilde x,t)=\tilde y+\Gamma^{(y,N)}(\tilde y,\tilde
x,t)\nonumber\\
u&=&u(\tilde y,\tilde x,t)=\tilde u+\Gamma^{(u,N)}(\tilde y,\tilde
x,t)\ .
\eeqa
The dissipative transformation $\Xi^{(N)}_d$ is defined by introducing suitable
functions with zero average over $\tilde x$ and $t$, say $\alpha^{(N)}$, $\beta^{(N)}$, $\gamma^{(N)}$ defined through series by
the coefficients $\alpha_{ji}$, $\beta_{ji}$, $\gamma_{ji}$, $j,i\in{\Z}_+$, such that
\beqa{B} X&=&\tilde x+\sum_{i=0}^N\sum_{j=0}^i \alpha_{j,i-j}
(\tilde y,\tilde x,t)\varepsilon^j\mu^{i-j}\equiv\tilde x+\alpha^{(N)}(\tilde y,\tilde x,t)\nonumber\\
Y&=&\tilde y+\sum_{i=0}^N\sum_{j=0}^i \beta_{j,i-j} (\tilde
y,\tilde x,t)\varepsilon^j\mu^{i-j}\equiv\tilde y+\beta^{(N)}(\tilde y,\tilde x,t)\nonumber\\
U&=&\tilde u+\sum_{i=0}^N\sum_{j=0}^i \gamma_{j,i-j} (\tilde
y,\tilde x,t)\varepsilon^j\mu^{i-j}\equiv\tilde u+\gamma^{(N)}(\tilde y,\tilde x,t)\ ,
\eeqa
with
$\alpha_{i0}(\tilde y,\tilde x,t)=\beta_{i0}(\tilde y,\tilde
x,t)=\gamma_{i0}(\tilde y,\tilde x,t)=0$ for any $i\geq 0$.
An iterative explicit construction of the vector functions $\psi_{j0}$, $\alpha_{j i}$, $\beta_{j i}$, $\gamma_{j i}$ will be given within
the proof of the Resonant Normal Form Lemma stated below.

\vskip.1in

\noindent \bf Remark 2. \rm The normal form equation defining the
generating function $\psi_{j0}(\tilde y,\tilde x,t)$ at order $j$
is given by
$$
\omega(\tilde y)\ \psi_{j0,x}(\tilde y,\tilde x,t)
+\psi_{j0,t}(\tilde y,\tilde x,t)+L_{j0}^{(nr,\leq K)}(\tilde y,\tilde x,t)=0\ ,
$$
for a suitable known function $L_{j0}^{(nr,\leq K)}(\tilde y,\tilde x,t)$ with zero average over $(\tilde x,t)$
and not containing resonant terms, say
$$
L_{j0}^{(nr,\leq K)}(\tilde y,\tilde x,t)=\sum_{(k,j)\in{\Z}^{\ell+1}\backslash\Lambda,\, |k|+|j|\leq K} \hat L_{j0,kj}(\tilde y)
\, e^{i(k\cdot\tilde x+jt)}\ .
$$
This equation can be solved provided $\omega=\omega(\tilde y)$ satisfies a
non--resonance condition of the form
$$
\omega(\tilde y)\cdot k+j\not=0\qquad {\rm for\ all}\ (k,j)\in{\Z}^{\ell+1}\backslash\Lambda\ ,\
|k|+|j|\leq K\ ,
$$
which is guaranteed by \equ{omega0}, provided $\varepsilon$
satisfies a smallness condition. Analogously, the dissipative
normal form provides an explicit construction of the functions
$\alpha^{(N)}$, $\beta^{(N)}$, $\gamma^{(N)}$, thanks to a
suitable choice of the drifts $\eta$, $\sigma$ and to the
assumption \equ{omega0}. More precisely, once expressed in terms
of the new variables $(Y,X,t)$, the functions $\beta_{j i}$ must
satisfy a normal form equation of the form \beq{beta} \omega(Y)
\beta_{j i,x}(Y,X,t)+\beta_{j i,t}(Y,X,t)+N_{j i}^{(nr,\leq
K)}(Y,X,t)+ \bar N_{j i}(Y)+N_{j i}^{(r,\leq
K)}(Y,X,t)+\eta_{ji}(Y,X,t)=0\ , \eeq for some known function
$N_{j i}\equiv \bar N_{j i}+ N_{j i}^{(nr,\leq K)}+ N_{j
i}^{(r,\leq K)}+N_{j i}^{(>K)}$; therefore, equation \equ{beta}
can be solved provided the drift components $\eta_{ji}(Y,X,t)$ are
chosen as the opposite of the sum of the average and of the
resonant parts:
$$
\eta_{ji}(Y,X,t)=-\Big(\bar N_{j i}(Y)+N_{j i}^{(r,\leq K)}(Y,X,t)\Big)\ .
$$
An analogous relation holds for $\gamma_{ji}$ and $\sigma_{ji}$.
This explains why the drift must be properly defined in order to be able to build the coordinate transformation
\equ{c}. This is not unusual, but it happens also in KAM theory (see e.g.\cite{CCL}).
We proceed now to state the Resonant Normal Form Lemma, which extends the Normal Form Lemma of
\cite{cellho2010a} to the resonant case of a frequency vector satisfying \equ{res}, \equ{omega0}.

\vskip.2in

\noindent \bf Resonant Normal Form Lemma. \sl Consider the vector
field \equ{SE} analytic in the complex extension
$C_{r_0}(A)\times C_{s_0}({\T}^{\ell+1})$ for some $r_0$, $s_0>0$.
Consider the extended vector field \equ{SEexp} on $A\times{\R}\times{{\T}^{\ell+1}}$.
For a given lattice $\Lambda\subset{\Z}^{\ell+1}$,
let $y_0\in A$, $K\in{\Z}_+$, $D\subseteq A$, $a>0$ be such that \equ{res} and \equ{omega0} are satisfied.
There exist suitable drift functions $\eta=\eta(y,x,t)$, $\sigma=\sigma(y,x,t)$, and there exist
$\varepsilon_0$, $\mu_0>0$ depending on $r_0$, $s_0$, $K$, $a$ and on the norms of
$\omega$, $h$, $f$, $g$, such that for
any $(\varepsilon,\mu)\leq(\varepsilon_0,\mu_0)$, one can construct a change of variables
close to the identity, say $\Xi^{(N)}:A\times{\R}\times{\T}^{\ell+1}\to{\R}^{\ell+1}\times{\T}^{\ell+1}$
with $(Y,X,U,t)=\Xi^{(N)}(y,x,u,t)$, being $(Y,U)\in{\R}^{\ell+1}$, $(X,t)\in{\T}^{\ell+1}$,
$N\in\Z_+$, which transforms \equ{SEexp} into a normal
form of order $N$ as in \equ{NF}.
Let $R_0<r_0$, $S_0<s_0$; having set $\lambda=\max(\eps,\mu)$, the normalized frequency is bounded by
\beq{omegadN}
\|\Omega_d^{(N)}-\omega\|_{R_0}\leq C_\omega\lambda\ ,
\eeq
where $C_\omega$ is a positive constant depending on $r_0$, $N$ and on the norms of $\omega$,
$h$, $f$, $g$. Denoting by $\Pi_y$ the projection on the $y$--coordinate, one gets
\beq{Pi} \|\Pi_y(\Xi_d^{(N)}\circ\Xi_c^{(N)})-Id\|\leq C_p\lambda\ , \eeq for some
positive constant $C_p$ depending on $r_0$, $s_0$, $N$ and on the norms of $\omega$, $h$, $f$, $g$.
With reference to the normal form
\equ{NF}, one has the following estimate
\beq{FG}
\|G^{(r,\leq K)}\|_{R_0,S_0}+\|G_{N+1}\|_{R_0,S_0}+\|G^{(>K)}\|_{R_0,S_0}\leq \lambda\,G+C_Y\lambda^{N+1}\ ,
\eeq
for some constant $C_Y$ and having bounded $\|G^{(r,\leq K)}\|_{R_0,S_0}$ by $\lambda\,G$, where
$G$ and $C_Y$ depend on $r_0$, $s_0$, $N$, $K$ and on the norms of $\omega$, $h$, $f$, $g$.
Choosing\footnote{The choice of $N$ is motivated as follows.
The relation $\lambda^N=e^{-K\tau_0}$ implies $N\log
\lambda=-K\tau_0$, namely $N=[K\tau_0/|\log \lambda|]$, where
$[\cdot]$ denotes the integer part.} $N=[K\tau_0/|\log \lambda|]$
for some $\tau_0>0$, one obtains that \equ{FG} becomes
\beq{18bis}
\|G^{(r,\leq K)}\|_{R_0,S_0}+\|G_{N+1}\|_{R_0,S_0}+\|G^{(>K)}\|_{R_0,S_0}
\leq \lambda\,G+C_Y\lambda e^{-K\tau_0}\ .
\eeq
\rm

\vskip.2in

\noindent
Before giving the proof of the Lemma, we provide the
statement of the main result, namely a bound on the variation of the
variables which are actions of the conservative system.
The following Theorem will be obtained through the
Resonant Normal Form Lemma under the resonance condition \equ{res}
and the quasi--convexity assumption \equ{QC}. Let us write the normal
form equations \equ{NF} using the following notation:
\beqa{NFcompact0}
\dot X &=& \omega(Y)+\varepsilon p_Y^{(\leq K)}(Y,X,t)+\mu s^{(\leq K)}(Y,X,t)
+F_{N+1}(Y,X,t)+F^{(>K)}(Y,X,t)\nonumber\\
\dot Y &=& -\varepsilon p_X^{(\leq K)}(Y,X,t)+G_{N+1}(Y,X,t)+G^{(>K)}(Y,X,t)\nonumber\\
\dot U &=& -\varepsilon p_t^{(\leq K)}(Y,X,t)+H_{N+1}(Y,X,t)+H^{(>K)}(Y,X,t)\ ,
\eeqa
where $p_X^{(\leq K)}$, $p_Y^{(\leq K)}$, $p_t^{(\leq K)}$ (independent of $\mu$) are the resonant contributions stemming just
from the conservative transformation, while $s^{(\leq K)}$ (depending on $\mu$ and $\varepsilon$)
represents the resonant part coming from the dissipative transformation.

\vskip.2in

\noindent
\bf Theorem \sl Consider the vector field \equ{SE} defined on $A\times {\T}^{\ell+1}$, satisfying the quasi--convexity
assumption \equ{QC}.
Let $y_0\in A$, $K\in{\Z}_+$, $D\subseteq A$, $a>0$ be such that \equ{res} and \equ{omega0} are satisfied.
Assume there exists $\varepsilon_0$, $\mu_0$, such that for
$(\varepsilon,\mu)\leq(\varepsilon_0,\mu_0)$, the Resonant Normal Form Lemma holds.
Let $\tau_0$, $C_p$, $\lambda$, $r_0$, $s_0$ as in the Resonant Normal Form Lemma.
With reference to \equ{NFcompact0}, we have that:
\begin{description}
    \item $i)$ if $p_X^{(\leq K)}=0$ or $s^{(\leq K)}=0$, then there exist $\rho_1>0$, $C_0>0$, such that
    $\|y(t)-y(0)\|\leq 2C_p\lambda+\rho_1$ for $t\leq T_1\equiv C_0 e^{K\tau_0}$, where $C_0$ depends on
$M$, $r_0$, $s_0$, $K$, $N$ and on the norms of $\omega$, $h$, $f$, $g$, while $\rho_1$ depends on the
above and on $m$, $\lambda$, $\Lambda$;
    \item $ii)$ if $p_X^{(\leq K)}\not=0$ and $s^{(\leq K)}\not=0$, then there exist $\rho_2>0$, $C_0'>0$, $C_0''>0$, such that if
    $t\leq T_2\equiv\min(C_0'e^{K\tau_0}, {{C_0''}\over {\varepsilon\mu}})$, then
    $\|y(t)-y(0)\|\leq 2C_p\lambda+\rho_2$, where $C_0'$, $C_0''$ depend on $M$, $r_0$, $s_0$, $K$, $N$ and
on the norms of $\omega$, $h$, $f$, $g$, while $\rho_2$ depends on the above and on $m$, $\lambda$,
$\Lambda$.
   \end{description}
\rm

\vskip.2in

\noindent
\bf Remark 3. \rm The above theorem is stated in terms of the functions $p_X^{(\leq K)}$ and $s^{(\leq K)}$ appearing in
the normal form equations \equ{NFcompact0}; in order to decide which of the conditions $i)$ or $ii)$ of
the Theorem
is satisfied, one needs to know the explicit expression of the functions $f_{01}$, $g_{01}$, $h_{10}$
appearing in the vector field \equ{SE}, tracing the resonant terms which could generate $p_X^{(\leq K)}$ and $s^{(\leq K)}$ by means of
an explicit construction of the normal form  or by means of
a \sl tree algorithm \rm (see, e.g., \cite{CFa}, \cite{GG2}, \cite{GJ}
and references therein).

\vskip.2in

\noindent
\bf Remark 4. \rm The Theorem states that in the non--resonant case (compare with \cite{cellho2010a}), as well as whenever the dissipative contribution to the resonant normal form is zero (at least up to the normalization
order), one finds a variation of the actions on exponential times; otherwise, there appears a fast drift of the actions on linear
(in $\varepsilon\,\mu$) times.

\vskip.1in

\noindent \bf Proof of the Resonant Normal Form Lemma. \rm By induction on the normalization order we prove that we can construct
a normal form of type \equ{NF} by means of suitable
transformations as in \equ{A} and \equ{B}. First we prove the statement by
constructing the first order normal form using the
conservative and then the dissipative transformation; next, we proceed to
construct the conservative and dissipative transformations at the order $N$. For sake of clarity,
we split the proof into four separate steps, referring, respectively, to the first order conservative and dissipative
normal forms, and to the $N$--th order conservative and dissipative transformations. Since the
conservative transformation is standard, we omit some details.

\vskip.1in

\noindent
\bf Step 1: Conservative transformation for $N=1$. \rm

\noindent
We start by implementing the first order transformation
\beqa{deltac1}
\tilde x&=&x+ \varepsilon \psi_{10,y}(\tilde y,x,t)\nonumber\\
y&=&\tilde y+ \varepsilon \psi_{10,x}(\tilde y,x,t)\nonumber\\
u&=&\tilde u+ \varepsilon \psi_{10,t}(\tilde y,x,t)\ ,
\eeqa
where $\psi_{10}=\psi_{10}(\tilde y,x,t)$ is an unknown function.
Let $\tilde r_0<r_0$, $\delta_0<s_0$, $\tilde s_0\equiv s_0-\delta_0$; then we can invert \equ{deltac1} as
\beqano
x&=&\tilde x+\varepsilon \Gamma^{(x,1)}(\tilde y,\tilde x,t)\nonumber\\
y&=&\tilde y+\varepsilon \Gamma^{(y,1)}(\tilde y,\tilde x,t)\nonumber\\
u&=&\tilde u+\varepsilon \Gamma^{(u,1)}(\tilde y,\tilde x,t)\ ,
\eeqano
for suitable functions $\Gamma^{(x,1)}$, $\Gamma^{(y,1)}$ and $\Gamma^{(u,1)}$, provided the following smallness condition on the
parameters is satisfied (compare with Appendix A):
\beq{C1}
70\ \|\psi_{10,y}\|_{\tilde r_0,s_0}\ e^{2s_0}\delta_0^{-1}\varepsilon<1\ .
\eeq
Using \equ{deltac1} and \equ{SEexp}, we obtain that the conservative normal form is achieved
whenever one can determine $\psi_{10}(\tilde y,\tilde x,t)$ such that
\beqa{otto}
\omega_y(\tilde y)\psi_{10,x}(\tilde y,\tilde x,t)+
\omega(\tilde y)\psi_{10,yx}(\tilde y,\tilde x,t)+\psi_{10,yt}(\tilde y,\tilde x,t)+ h^{(nr,\leq K)}_{10,y}(\tilde y,\tilde x,t)&=&0\nonumber\\
\omega(\tilde y)\psi_{10,xx}(\tilde y,\tilde x,t)+\psi_{10,xt}(\tilde y,\tilde x,t)+
 h^{(nr,\leq K)}_{10,x}(\tilde y,\tilde x,t)&=&0\nonumber\\
\omega(\tilde y)\psi_{10,tx}(\tilde y,\tilde x,t)+\psi_{10,tt}(\tilde y,\tilde x,t)+
 h^{(nr,\leq K)}_{10,t}(\tilde y,\tilde x,t)&=&0\ .
\eeqa
Let us define
$$
\Omega_c^{(1)}(\tilde y)=\Omega_c^{(1)}(\tilde y;\varepsilon)\equiv \omega(\tilde y)+\varepsilon \bar h_{10,y}(\tilde y) \ ;
$$
Equations \equ{otto} are equivalent to take the derivatives with respect to $y$,
$x$ and $t$ of
$$
\omega(\tilde y)\psi_{10,x}(\tilde y,\tilde x,t)+\psi_{10,t}(\tilde y,\tilde x,t)+
h^{(nr,\leq K)}_{10}(\tilde y,\tilde x,t)=0\ .
$$
Expanding $\psi_{10}$ and
$h^{(nr,\leq K)}_{10}$ into Fourier series, one obtains that $\psi_{10}$ is
given by the expression (independent of $\tilde u$):
\beq{psi10}
\psi_{10}(\tilde y,\tilde x,t)=\imath \sum_{(k,j)\in{\Z}^{\ell+1}\backslash\Lambda,\, |k|+|j|\leq K}
{{\hat {h}^{(nr,\leq K)}_{10,kj}(\tilde y)}\over {\omega(\tilde y)\cdot k+j}}\
e^{\imath(k\cdot \tilde x+jt)}\ .
\eeq
This function is well defined, since the zero and small divisors are controlled as follows.
The second of \equ{deltac1} can be inverted as $\tilde y=y+\varepsilon R^{(y,1)}(y,x,t)$ for a suitable
function $R^{(y,1)}=R^{(y,1)}(y,x,t)$, provided that for $\tilde r_0'<r_0$ one has (see Appendix A)
\beq{C2}
70\,\varepsilon\, \|\psi_{10,x}\|_{\tilde r_0,s_0}{1\over {\tilde r_0-\tilde r_0'}}<1\ ,
\eeq
being $\varepsilon\|R^{(y,1)}\|_{\tilde r_0',s_0}\leq\|\psi_{10,x}\|_{\tilde r_0,s_0}$.
Then, the divisors
appearing in \equ{psi10} are bounded by
\beq{omt}
|\omega(\tilde y)\cdot k+j|\geq a-\varepsilon K\|R^{(y,1)}\|_{\tilde r_0',s_0} \|\omega_y\|_{r_0}>{a\over 2}\ ,
\eeq
provided
\beq{C2bis}
\varepsilon <{a\over {2K\|R^{(y,1)}\|_{\tilde r_0',s_0} \|\omega_y\|_{r_0}}}\ .
\eeq

\vskip.1in

\noindent
\bf Step 2: Dissipative transformation for $N=1$. \rm

\noindent
We define the first--order dissipative transformation as
\beqa{D1}
X&=&\tilde x+\alpha_{01}(\tilde y,\tilde x,t)\mu\nonumber\\
Y&=&\tilde y+\beta_{01}(\tilde y,\tilde x,t)\mu\nonumber\\
U&=&\tilde u+\gamma_{01}(\tilde y,\tilde x,t)\mu\ ,
\eeqa
for unknown functions $\alpha_{01}$, $\beta_{01}$ and $\gamma_{01}$. Let us start by inverting \equ{D1} as
\beqa{starstar}
\tilde x&=&X-\alpha_{01}(Y,X,t)\mu+O_2(\mu)=X+\Delta^{(x,1)}(Y,X,t)\mu\nonumber\\
\tilde y&=&Y-\beta_{01}(Y,X,t)\mu+O_2(\mu)=Y+\Delta^{(y,1)}(Y,X,t)\mu\nonumber\\
\tilde u&=&U-\gamma_{01}(Y,X,t)\mu+O_2(\mu)=U+\Delta^{(u,1)}(Y,X,t)\mu
\eeqa
for suitable functions $\Delta^{(x,1)}$, $\Delta^{(y,1)}$ and $\Delta^{(u,1)}$ provided the following smallness conditions
on the parameters are satisfied (see Appendix A):
\beqa{C4}
70\ \|\alpha_{01}\|_{\tilde r_0,\tilde s_0}\ e^{2\tilde s_0}\tilde\delta_0^{-1}\, \mu&<&1\nonumber\\
70\ (\|\beta_{01}\|_{\tilde r_0,\tilde s_0}+\|\beta_{01,x}\|_{\tilde r_0,\tilde s_0}\ \|\alpha_{01}\|_{\tilde r_0,\tilde s_0})\
{1\over {\tilde r_0-R_0}}\, \mu&<&1\nonumber\\
70\ (\|\gamma_{01}\|_{\tilde r_0,\tilde s_0}+\|\gamma_{01,x}\|_{\tilde r_0,\tilde s_0}\ \|\alpha_{01}\|_{\tilde r_0,\tilde s_0})
\ {1\over {\tilde r_0-R_0}}\, \mu&<&1\ ,
\eeqa
where $\tilde \delta_0\equiv \tilde s_0/2$, $R_0<\tilde r_0$ and
being $\|\Delta^{(x,1)}\|_{\tilde r_0,\tilde s_0-\tilde \delta_0}\leq\|\alpha_{01}\|_{\tilde r_0,\tilde s_0}$.
Through \equ{D1} and \equ{starstar} we can express $\dot X$, $\dot Y$ as a function of
$X$, $Y$; the normal form is obtained assuming
that $\alpha_{01}$, $\beta_{01}$ and $\eta_{01}$ satisfy
the following equations: \beqa{XX}
\omega(Y)\alpha_{01,x}(Y,X,t)&+&\alpha_{01,t}(Y,X,t)-\omega_y(Y)\beta_{01}(Y,X,t)
+f_{01}^{(nr,\leq K)}(Y,X,t)=0\nonumber\\
\omega(Y)\beta_{01,x}(Y,X,t)&+&\beta_{01,t}(Y,X,t)\nonumber\\
&-&g^{(nr,\leq K)}_{01}(Y,X,t)-\bar g_{01}(Y)
-g^{(r,\leq K)}_{01}(Y,X,t)+\eta_{01}(Y,X,t)=0\nonumber\\
\omega(Y)\gamma_{01,x}(Y,X,t)&+&\gamma_{01,t}(Y,X,t)+\sigma_{01}(Y,X,t)=0\ .
\eeqa
Since $\alpha_{ij}$, $\beta_{ij}$, $\gamma_{ij}$ have zero average and they do not contain resonant terms, the system of equations \equ{XX} can be solved, provided that we choose $\eta_{01}(Y,X,t)$
as
$$
\eta_{01}(Y,X,t)\equiv\bar g_{01}(Y)+g^{(r,\leq K)}_{01}(Y,X,t)
$$
and that we set $\sigma_{01}=0$ as well as $\gamma_{01}=0$. Setting
$\Omega_d^{(1)}=\omega(Y)+\varepsilon{\bar h}_{10,y}(Y)+\mu \bar f_{01}(Y)$,
the first order normal form can be written as
\beqano
\dot X&=&\Omega_d^{(1)}
+\varepsilon h^{(r,\leq K)}_{10,y}(Y,X,t)+\mu f^{(r,\leq K)}_{01}(Y,X,t)\nonumber\\
&+&\varepsilon h^{(>K)}_{10,y}(Y,X,t)+\mu f^{(>K)}_{01}(Y,X,t)+F_2(Y,X,t)\nonumber\\
\dot Y&=&-\varepsilon h^{(r,\leq K)}_{10,x}(Y,X,t)-\varepsilon h^{(>K)}_{10,x}(Y,X,t)-\mu
g^{(>K)}_{01}(Y,X,t)+G_2(Y,X,t)\nonumber\\
\dot U&=&-\varepsilon h^{(r,\leq K)}_{10,t}(Y,X,t)-\varepsilon h^{(>K)}_{10,t}(Y,X,t)+H_2(Y,X,t)\ ,
\eeqano
where $F_2$, $G_2$ are functions of order $O_2(\varepsilon,\mu)$ and $H_2$ is a function of order $O_2(\varepsilon)$. We remark that
the solution of \equ{XX} involves small divisors of the form $\omega(Y)\cdot k+j$ with $(k,j)\in{\Z}^{\ell+1}$ and
$|k|+|j|\leq K$.
Using the same argument as in \equ{omt}, the small divisors are bounded by $a/4$ provided that the following smallness condition holds
(compare with Appendix A):
\beq{C5}
\mu <{a\over {4K\|\beta_{01}\|_{\tilde r_0,\tilde s_0} \|\omega_y\|_{r_0}}}\ .
\eeq

\vskip.1in

\noindent
\bf Step 3: Conservative transformation for the order $N$. \rm

\noindent
Assume that the Lemma holds to the order $N-1$. We introduce
the conservative transformation to the order $N$ as
\beqa{deltac2}
\tilde x&=&x+\sum_{j=1}^{N-1} \psi_{j0,y}(\tilde y,x,t)\varepsilon^j+\psi_{N0,y}(\tilde y,x,t)
\varepsilon^N\equiv x+\psi_y^{(N)}(\tilde y,x,t)\nonumber\\
y&=&\tilde y+\sum_{j=1}^{N-1} \psi_{j0,x}(\tilde y,x,t)\varepsilon^j+\psi_{N0,x}(\tilde y,x,t)
\varepsilon^N\equiv \tilde y+\psi_x^{(N)}(\tilde y,x,t)\nonumber\\
u&=&\tilde u+\sum_{j=1}^{N-1} \psi_{j0,t}(\tilde y,x,t)\varepsilon^j+\psi_{N0,t}(\tilde y,x,t)
\varepsilon^N\equiv \tilde u+\psi_t^{(N)}(\tilde y,x,t)\ ,
\eeqa
where for $1\leq j\leq N-1$ the functions $\psi_{j0}(\tilde y,x,t)$ are assumed to be known.
We can invert \equ{deltac2} as
\beqa{deltac3}
x&=&x(\tilde y,\tilde x,t)\nonumber\\
&=&\tilde x+\sum_{j=1}^{N} \Gamma_{j0}^{(x)}(\tilde y,\tilde x,t)\varepsilon^j-
\psi_{N0,y}(\tilde y,\tilde x,t)\varepsilon^N+O_{N+1}(\varepsilon)
\equiv \tilde x+\Gamma^{(x,N)}(\tilde y,\tilde x,t)\nonumber\\
y&=&y(\tilde y,\tilde x,t)\nonumber\\
&=&\tilde y+\sum_{j=1}^{N} \Gamma_{j0}^{(y)}(\tilde y,\tilde x,t)\varepsilon^j+
\psi_{N0,x}(\tilde y,\tilde x,t)\varepsilon^N+O_{N+1}(\varepsilon)\equiv \tilde y+\Gamma^{(y,N)}(\tilde y,\tilde x,t)\nonumber\\
u&=&u(\tilde y,\tilde x,t)\nonumber\\
&=&\tilde u+\sum_{j=1}^{N} \Gamma_{j0}^{(u)}(\tilde y,\tilde x,t)\varepsilon^j+
\psi_{N0,t}(\tilde y,\tilde x,t)\varepsilon^N+O_{N+1}(\varepsilon)\equiv \tilde u+\Gamma^{(u,N)}(\tilde y,\tilde x,t)\ ,
\eeqa
provided that (see Appendix A), choosing
$\tilde r_0<r_0$, $\delta_0<s_0$:
\beq{33ter}
70\|\psi^{(N)}_y\|_{\tilde r_0,s_0}\ e^{2s_0}\delta_0^{-1}\ <\ 1\ ,
\eeq
being $\Gamma_{j0}^{(x)}$, $\Gamma_{j0}^{(y)}$, $\Gamma_{j0}^{(u)}$, $1\leq j\leq N$,
known functions.
We proceed to compute $\dot x$, $\dot y$, $\dot u$ as a function of
$\tilde x$, $\tilde y$, $t$ and after expanding in Taylor series, we obtain
\beqa{AA}
\dot x&=&\omega(\tilde y)+\omega_y(\tilde y)\psi_{N0,x}(\tilde y,\tilde x,t)\varepsilon^N+
F^{(1,\leq K,\leq N)}(\tilde y,\tilde x,t)+O_{N+1}^{>K}\nonumber\\
\dot y&=&G^{(1,\leq K,N)}(\tilde y,\tilde x,t)+\mu\left(\eta_{N-1,1}(\tilde y,\tilde x,t)\varepsilon^{N-1}+...+
\eta_{0,N}(\tilde y,\tilde x,t)\mu^{N-1}\right)+O_{N+1}^{>K}\nonumber\\
\dot u&=&H^{(1,\leq K,N)}(\tilde y,\tilde x,t)+\mu\left(\sigma_{N-1,1}(\tilde y,\tilde x,t)\varepsilon^{N-1}+...+
\sigma_{0,N}(\tilde y,\tilde x,t)\mu^{N-1}\right)+O_{N+1}^{>K}\ ,
\eeqa
where $O_{N+1}^{>K}$ is a compact notation to denote terms of order $O_{N+1}(\varepsilon,\mu)$ and/or containing Fourier components
greater than $K$; the functions $F^{(1,\leq K,\leq N)}$, $G^{(1,\leq K,N)}$, $H^{(1,\leq K,N)}$ are known,
contain Fourier components up to the order $K$, contain orders in $\varepsilon$ and $\mu$ up to the order
$N$ and they are at most linear in $\mu$. Using \equ{deltac2}, \equ{AA} and the inductive hypothesis,
the conservative normal form at the order $N$ is obtained once the function
$\psi_{N0}$ satisfies the equations
\beqa{psi}
\omega_y(\tilde y)\psi_{N0,x}(\tilde y,\tilde x,t)+\omega(\tilde y)
\psi_{N0,yx}(\tilde y,\tilde x,t)+\psi_{N0,yt}(\tilde y,\tilde x,t)+
F^{(2,nr,\leq K,\leq N)}_{N0}(\tilde y,\tilde x,t)&=&0\nonumber\\
\omega(\tilde y)\psi_{N0,xx}(\tilde y,\tilde x,t)+\psi_{N0,xt}(\tilde y,\tilde x,t)-
G^{(2,nr,\leq K,\leq N)}_{N0}(\tilde y,\tilde x,t)&=&0\nonumber\\
\omega(\tilde y)\psi_{N0,tx}(\tilde y,\tilde x,t)+\psi_{N0,tt}(\tilde y,\tilde x,t)-
H^{(2,nr,\leq K,\leq N)}_{N0}(\tilde y,\tilde x,t)&=&0\ ,\nonumber\\
\eeqa
where $F^{(2,nr,\leq K,\leq N)}_{N0}$, $G^{(2,nr,\leq K,\leq N)}_{N0}$, $H^{(2,nr,\leq K,\leq N)}_{N0}$
are the non--resonant parts of known functions $F^{(2,\leq K, \leq N)}$,  $G^{(2,\leq K, \leq N)}$
 $H^{(2,\leq K, \leq N)}$ that we decompose as
\beqa{starf2}
F^{(2,\leq K,\leq N)}(\tilde y,\tilde x,t)&\equiv& \sum_{j=1}^N \bar F^{(2,\leq N)}_{j0}(\tilde y)\varepsilon^j+
\sum_{j=1}^{N} F^{(2,nr,\leq K,\leq N)}_{j0}(\tilde y,\tilde x,t)\varepsilon^j\nonumber\\
&+&\sum_{j=1}^{N} F^{(2,r,\leq K,\leq N)}_{j0}(\tilde y,\tilde x,t)\varepsilon^j
+\mu \sum_{j=0}^{N-1} F^{(2,\leq K,\leq N)}_{j1}(\tilde y,\tilde x,t)\varepsilon^j\
\eeqa
(and similar for the remaining functions). From the Hamiltonian structure it can be easily recognized that $F^{(2,nr,\leq K,\leq N)}_{N0}$, $G^{(2,nr,\leq K,\leq N)}_{N0}$, $H^{(2,nr,\leq K,\leq N)}_{N0}$ are,
respectively, the derivatives with respect to $y$, $x$, $t$ of the same function, so that equations
\equ{psi} uniquely define the solution $\psi_{N0}(\tilde y,\tilde x,t)$.
We are finally led to the following conservative normal form:
\beqa{CNF}
\dot{\tilde x}&=&\Omega_c^{(N)}(\tilde y;\varepsilon)+\sum_{j=1}^N F_{j0}^{(2,r,\leq K,\leq N)}(\tilde y,\tilde x,t)\varepsilon^j
+\mu \sum_{j=0}^{N-1} F_{j1}^{(2,\leq K,\leq N)}(\tilde y,\tilde x,t)\varepsilon^j
+O_{N+1}^{>K}\nonumber\\
\dot{\tilde y}&=&\sum_{j=1}^N G_{j0}^{(2,r,\leq K,\leq N)}(\tilde y,\tilde x,t)\varepsilon^j+
\mu \sum_{j=0}^{N-1} G_{j1}^{(2,\leq K,\leq N)}(\tilde y,\tilde x,t)\varepsilon^j
+O_{N+1}^{>K}\nonumber\\
\dot{\tilde u}&=&\sum_{j=1}^N H_{j0}^{(2,r,\leq K,\leq N)}(\tilde y,\tilde x,t)\varepsilon^j+
\mu \sum_{j=0}^{N-1} H_{j1}^{(2,\leq K,\leq N)}(\tilde y,\tilde x,t)\varepsilon^j
+O_{N+1}^{>K}\ ,
\eeqa
where
$$
\Omega_c^{(N)}(\tilde y;\varepsilon)\equiv\omega(\tilde y)+\sum_{j=1}^N\bar F_{j0}^{(2,N)}(\tilde y)
\varepsilon^j\ ,
$$
which implies that $\|\Omega_c^{(N)}-\omega\|\leq C_c\varepsilon$
for a suitable constant $C_c$. The normal form equations can be
solved, provided that the small divisors taking the expression
$\omega(\tilde y)\cdot k+j$, for $k\in{\Z}^\ell$, $j\in{\Z}$ with
$|k|+|j|\leq K$, are controlled by a non--resonance condition,
which is guaranteed whenever (see Appendix A) \beq{C6}
\varepsilon\leq{a\over {2K\|R^{(y,N)}\|_{\tilde
r_0',s_0}\|\omega_y\|_{r_0}}}\ , \eeq where $R^{(y,N)}$ is the
function inverting the transformation, namely $\tilde
y=y+\varepsilon R^{(y,N)}(y,x,t)$, and $\tilde r_0'<\tilde r_0$.
The inversion can be performed provided (see Appendix A)
\beq{cnew1} 70\, \|\psi_x^{(N)}\|_{\tilde r_0,s_0}\ {1\over
{r_0-\tilde r_0'}}\ <\ 1 \eeq with $\|R^{(y,N)}\|_{\tilde
r_0',s_0}\leq \|\psi_x^{(N)}\|_{\tilde r_0,s_0}$.

\vskip.1in

\noindent
\bf Step 4: Dissipative transformation for the order $N$. \rm

\noindent
We consider the transformation \equ{B} at the order $N$, which can be inverted as
\beqa{star5}
\tilde x&=&\tilde x(Y,X,t)\nonumber\\
&=&X+\sum_{i=0}^{N-2}\sum_{j=1}^{N-1-i} a_{ij}(Y,X,t)\varepsilon^i\mu^j-
\sum_{i=0}^{N-1} \alpha_{i,N-i}(Y,X,t)\varepsilon^i\mu^{N-i}+O_{N+1}(\varepsilon,\mu)\nonumber\\
\tilde y&=&\tilde y(Y,X,t)\nonumber\\
&=&Y+\sum_{i=0}^{N-2}\sum_{j=1}^{N-1-i} b_{ij}(Y,X,t)\varepsilon^i\mu^j-
\sum_{i=0}^{N-1} \beta_{i,N-i}(Y,X,t)\varepsilon^i\mu^{N-i}+O_{N+1}(\varepsilon,\mu)\nonumber\\
\tilde u&=&\tilde u(Y,X,t)\nonumber\\
&=&U+\sum_{i=0}^{N-2}\sum_{j=1}^{N-1-i} c_{ij}(Y,X,t)\varepsilon^i\mu^j-
\sum_{i=0}^{N-1} \gamma_{i,N-i}(Y,X,t)\varepsilon^i\mu^{N-i}+O_{N+1}(\varepsilon,\mu)\ ,
\eeqa
for suitable known functions $a_{ij}(Y,X,t)$, $b_{ij}(Y,X,t)$, $c_{ij}(Y,X,t)$,
provided that the parameters satisfy (see Appendix A):
\beqa{C7}
70\|\alpha^{(N)}\|_{\tilde r_0,\tilde s_0}\ e^{2\tilde s_0}\tilde \delta_0^{-1}&<&1\nonumber\\
70\left(\|\beta^{(N)}\|_{\tilde r_0,\tilde s_0}+\|\beta_x^{(N)}\|_{\tilde r_0,\tilde s_0}\|\alpha^{(N)}\|_{\tilde r_0,\tilde s_0}\right)\ {1\over {\tilde r_0-R_0}}&<&1\nonumber\\
70\left(\|\gamma^{(N)}\|_{\tilde r_0,\tilde s_0}+\|\gamma_x^{(N)}\|_{\tilde r_0,\tilde s_0}\|\alpha^{(N)}\|_{\tilde r_0,\tilde s_0}\right)\ {1\over {\tilde r_0-R_0}}&<&1\ ,
\eeqa
where $\tilde\delta_0\equiv \tilde s_0/2$, $R_0<\tilde r_0$.
In order to determine the unknown
functions $\alpha_{0,N}$, ..., $\alpha_{N-1,1}$, $\beta_{0,N}$, ..., $\beta_{N-1,1}$, $\gamma_{0,N}$, ..., $\gamma_{N-1,1}$,
$\eta_{N-1,0}$, ..., $\eta_{0,N-1}$, $\sigma_{N-1,0}$, ..., $\sigma_{0,N-1}$, using \equ{CNF} and \equ{star5} we express
$\dot{\tilde x}$, $\dot{\tilde y}$ in terms of $X$, $Y$ and we compute $\dot X$, $\dot Y$ using \equ{B}, \equ{CNF}, \equ{star5} as
\beqano
\dot X&=&\omega(Y)-\omega_y(Y)\left(\sum_{i=0}^{N-1}\beta_{i,N-i}(Y,X,t)\varepsilon^i\mu^{N-i}\right)
+\omega(Y)\sum_{i=0}^{N-1}\alpha_{i,N-i,x}(Y,X,t)\varepsilon^i\mu^{N-i}\nonumber\\
&+&\sum_{i=0}^{N-1}\alpha_{i,N-i,t}(Y,X,t)\varepsilon^i\mu^{N-i}+F^{(3,nr,\leq K,N)}(Y,X,t)+\bar F^{(3,\leq N)}(Y)\nonumber\\
&+&F^{(3,r,\leq K,\leq N)}(Y,X,t)+O_{N+1}^{>K}\nonumber\\
\dot Y&=&\mu\Big(\eta_{N-1,1}(Y,X,t)\varepsilon^{N-1}+ ...+\eta_{0,N}(Y,X,t)\mu^{N-1}\Big)
+\omega(Y)\sum_{i=0}^{N-1}\beta_{i,N-i,x}(Y,X,t)\varepsilon^i\mu^{N-i}\nonumber\\
&+&\sum_{i=0}^{N-1}\beta_{i,N-i,t}(Y,X,t)\varepsilon^i\mu^{N-i}+G^{(3,nr,\leq K,N)}(Y,X,t)+\bar G^{(3,N)}(Y)\nonumber\\
&+&\sum_{i=1}^N G_{i0}^{(3,r,\leq K,\leq N)}(Y,X,t)\varepsilon^i
+G^{(3,r,\leq K,N)}(Y,X,t)+O_{N+1}^{>K}\nonumber\\
\dot U&=&\mu\Big(\sigma_{N-1,1}(Y,X,t)\varepsilon^{N-1}+ ...+\sigma_{0,N}(Y,X,t)\mu^{N-1}\Big)
+\omega(Y)\sum_{i=0}^{N-1}\gamma_{i,N-i,x}(Y,X,t)\varepsilon^i\mu^{N-i}\nonumber\\
&+&\sum_{i=0}^{N-1}\gamma_{i,N-i,t}(Y,X,t)\varepsilon^i\mu^{N-i}+H^{(3,nr,\leq K,N)}(Y,X,t)+\bar H^{(3,N)}(Y)\nonumber\\
&+&\sum_{i=1}^N H_{i0}^{(3,r,\leq K,\leq N)}(Y,X,t)\varepsilon^i+H^{(3,r,\leq K,N)}(Y,X,t)+O_{N+1}^{>K}\ ,
\eeqano
where $F^{(3,nr,\leq K,N)}$, $G^{(3,nr,\leq K,N)}$, $H^{(3,nr,\leq K,N)}$ denote known non--resonant functions that
we can expand as
$$
F^{(3,nr,\leq K,N)}(Y,X,t)=\sum_{i=0}^{N-1} F^{(3,nr,\leq K,N)}_{i,N-1}(Y,X,t)\varepsilon^i\mu^{N-i}
$$
and similarly for $G^{(3,nr,\leq K,N)}$, $H^{(3,nr,\leq K,N)}$;
$F^{(3,r,\leq K,\leq N)}$, $G^{(3,r,\leq K,N)}$, $H^{(3,r,\leq
K,N)}$ denote known resonant functions; $\bar F^{(3,\leq N)}$,
$\bar G^{(3,N)}$, $\bar H^{(3,N)}$ denote the average terms.
Recall that due to the inductive hypothesis, the functions
$\alpha_{ij}$, $\beta_{ij}$, $\gamma_{ij}$, $\eta_{ij}$,
$\sigma_{ij}$, determine a normal form up to the order
$\varepsilon^i\mu^j$ with $0\leq i+j\leq N-1$. The normal form at
the order $N$ is obtained by imposing that $\alpha_{ij}$,
$\beta_{ij}$, $\gamma_{ij}$, $\eta_{ij}$, $\sigma_{ij}$ satisfy
the normal form equations
\beqa{star7}
&&-\omega_y(Y)\beta_{i,N-i}(Y,X,t)+\omega(Y)\alpha_{i,N-i,x}(Y,X,t)+
\alpha_{i,N-i,t}(Y,X,t)\nonumber\\
&&\qquad\qquad\qquad\qquad\qquad\qquad\qquad\qquad\qquad\qquad+F^{(3,nr,\leq K,N)}_{i,N-i}(Y,X,t)=0\nonumber\\
&&\omega(Y)\beta_{i,N-i,x}(Y,X,t)+
\beta_{i,N-i,t}(Y,X,t)+G_{i,N-i}^{(3,nr,\leq K,N)}(Y,X,t)
+G_{i,N-i}^{(3,r,\leq K,\leq N)}(Y,X,t)\nonumber\\
&&\qquad\qquad\qquad\qquad\qquad\qquad\qquad\qquad\qquad+\bar G_{i,N-i}^{(3,N)}(Y)+\eta_{i,N-i}(Y,X,t)=0
\nonumber\\
&&\omega(Y)\gamma_{i,N-i,x}(Y,X,t)+
\gamma_{i,N-i,t}(Y,X,t)+H_{i,N-i}^{(3,nr,\leq K,N)}(Y,X,t)
+H_{i,N-i}^{(3,r,\leq K,\leq N)}(Y,X,t)\nonumber\\
&&\qquad\qquad\qquad\qquad\qquad\qquad\qquad\qquad\qquad+\bar H_{i,N-i}^{(3,N)}(Y)+\sigma_{i,N-i}(Y,X,t)=0
\eeqa
for $0\leq i\leq N-1$. Equations \equ{star7} can be solved provided $\omega(Y)$ satisfies a
non--resonance condition, which is guaranteed by (see Appendix A)
\beq{C8}
4K\|\omega_y\|_{r_0} \|\ \beta^{(N)}\|_{\tilde r_0,\tilde s_0} <a\ ,
\eeq
where we intend that $Y\equiv \tilde y+\beta^{(N)}(\tilde y,\tilde x,t;\varepsilon,\mu)$. From the second and third
of \equ{star7}, we get
\beqano
\eta_{i,N-i}(Y,X,t)&\equiv& \bar G^{(3,N)}_{i,N-i}(Y)+G_{i,N-i}^{(3,r,\leq K,N)}(Y,X,t)\varepsilon^i\mu^{N-i}\nonumber\\
\sigma_{i,N-i}(Y,X,t)&\equiv& \bar H^{(3,N)}_{i,N-i}(Y)+H_{i,N-i}^{(3,r,\leq K,N)}(Y,X,t)\varepsilon^i\mu^{N-i}\ .
\eeqano
Setting
\beq{omegad}
\Omega_d^{(N)}(Y;\varepsilon,\mu)\equiv \omega(Y)+\sum_{i=1}^N\Omega_{i0}(Y)\varepsilon^i+
\sum_{j=1}^N\sum_{i=0}^{N-j}\Omega_{ij}(Y)\varepsilon^i\mu^j\ ,
\eeq
with $\Omega_{i0}\equiv \bar F_{i0}^{(3,\leq N)}(Y)$ and $\Omega_{i,N-i}\equiv \bar F^{(3,\leq N)}_{i,N-i}(Y)$,
the normal form is finally given by
\beqa{NFf}
\dot X &=& \Omega_d^{(N)}(Y;\varepsilon,\mu)+\sum_{i=1}^N F_{i0}^{(3,r,\leq K,\leq N)}(Y,X,t)\varepsilon^i
+\sum_{j=1}^N \sum_{i=0}^{N-j} F_{ij}^{(3,r,\leq K,\leq N)}(Y,X,t)\varepsilon^i\mu^j\nonumber\\
&+&F_{N+1}(Y,X,t)+F^{(>K)}(Y,X,t)\nonumber\\
\dot Y  &=& \sum_{i=1}^N G_{i0}^{(3,r,\leq K,\leq N)}(Y,X,t)\varepsilon^i+G_{N+1}(Y,X,t)+G^{(>K)}(Y,X,t)\nonumber\\
\dot U  &=& \sum_{i=1}^N H_{i0}^{(3,r,\leq K,\leq N)}(Y,X,t)\varepsilon^i+H_{N+1}(Y,X,t)+H^{(>K)}(Y,X,t)\ ,
\eeqa
where $F_{N+1}$, $G_{N+1}$ are $O_{N+1}(\varepsilon,\mu)$, $H_{N+1}$ is order $O_{N+1}(\varepsilon)$ and $F^{(>K)}$, $G^{(>K)}$, $H^{(>K)}$ contain only
terms with Fourier index greater than $K$. The normal form \equ{NF} is recovered  with an obvious
identification of the functions $F^{(r,\leq K)}$, $G^{(r,\leq K)}$, $H^{(r,\leq K)}$. The smallness
requirements on $\varepsilon$, $\mu$, say $\varepsilon\leq \varepsilon_0$,
$\mu\leq\mu_0$, are needed to guarantee the non--resonance condition
(see \equ{C2bis}, \equ{C5}, \equ{C6}, \equ{C8}) and the inversion of the transformations
(see \equ{C1}, \equ{C2}, \equ{C4}, \equ{33ter}, \equ{cnew1}, \equ{C7}). The estimate \equ{omegadN}
holds true, due to the definition of $\Omega_d^{(N)}$ in \equ{omegad}.
The estimate \equ{Pi} follows from the fact that \equ{ast1} is close to the identity up to first order.

\vskip .1in

\noindent
Due to the exponential decay of the Fourier coefficients
(compare with Lemma B.1 of Appendix B), we can bound $G^{(>K)}$ for some $\tau_0>0$ as
\beq{cgtilde}
\|G^{(>K)}\|_{R_0,S_0}\leq \tilde C_G \lambda e^{-K\tau_0}
\eeq
for a suitable constant $\tilde C_G$. On the other hand we can bound $G_{N+1}$ in \equ{NFf} as
\beq{cg}
\|G_{N+1}\|_{R_0,S_0}\leq C_G\lambda^{N+1}\ ,
\eeq
for a suitable constant $C_G$. Finally, from the second of \equ{NFf} we obtain:
\beqano
\|\sum_{i=1}^N G_{i0}^{(3,r,\leq K,\leq N)}(Y,X,t)\varepsilon^j\|_{R_0,S_0}+\|G_{N+1}\|_{R_0,S_0}&+&\|G^{(>K)}\|_{R_0,S_0}\nonumber\\
&\leq& \lambda\ G+\lambda^{N+1}C_G+\lambda\tilde C_G  e^{-K\tau_0}\ ,
\eeqano
having defined $\lambda G$ as an upper bound of $\sup_{\eps\leq\eps_0}\|\sum_{i=1}^N G_{i0}^{(3,r,\leq K,\leq N)}(Y,X,t)\varepsilon^i\|_{R_0,S_0}$.
Choosing $N=[K\tau_0/|\log\lambda|]$, we obtain \equ{FG} and \equ{18bis} with $C_Y\equiv C_G+\tilde C_G$.
This concludes the proof of the Lemma. $\Box$

\vskip.2in

\noindent
\bf Proof of the theorem. \rm The distance between $y(t)$ and the initial
condition $y(0)$ for $t\geq 0$ can be bounded by the sum of the following terms:
\beq{yY}
\|y(t)-y(0)\|\leq \|y(t)-Y(t)\|+\|Y(t)-Y(0)\|+\|Y(0)-y(0)\|\ .
\eeq
By the estimate \equ{Pi} of the Resonant Normal Form Lemma, one obtains
$$
\|y(t)-Y(t)\|\leq C_p\,\lambda\ ,\qquad \|y(0)-Y(0)\|\leq C_p\,\lambda\ .
$$
By the second of \equ{NF} and by \equ{FG}, one gets
\beqano
\|Y(t)-Y(0)\|&\leq& \int_0^t \Big(\|G^{(r,\leq K)}\|_{R_0,S_0}+\|G_{N+1}\|_{R_0,S_0}+\|G^{(>K)}\|_{R_0,S_0}\Big)ds\nonumber\\
&\leq& \lambda\,G\,t+C_Y\lambda^{N+1} t\ ,
\eeqano
which indicates that the action variation takes place on linear time scales due to the term $\lambda\,G\,t$, while exponential
times are associated to the term $C_Y\,\lambda^{N+1}\,t$.
We remark that $G=0$ corresponds to the absence of resonant terms in the normal form for $Y$.
Notice that the case of non--resonant stability estimates given in \cite{cellho2010a}
is recovered whenever also the resonant terms in the $X$ variable are zero.
Let us start with the case $G=0$; for a suitable $\rho_1>0$ that we write as $\rho_1=C_\rho \lambda$
for some $C_\rho>0$, let
$$
t\leq {C_\rho\over {C_Y}}\ e^{K\tau_0}\ .
$$
Finally, setting $\rho_0\equiv (2C_p+C_\rho)\,\lambda$, we obtain the following variation of the original action variables on exponential times:
$$
\|y(t)-y(0)\|\leq \rho_0\qquad {\rm for\ } t\leq C_0 e^{K\tau_0}\ ,
$$
having defined $C_0\equiv C_\rho/C_Y$. This result is in agreement with statement $i)$, once $G^{(r,\leq K)}$ is identified
with $-\varepsilon p_X^{(\leq K)}$ (compare with \equ{NFcompact0}).

\noindent
Next, we study the case $G\not=0$; to this end, we compute the variation of the \sl energy \rm (i.e. the Lyapunov function, see e.g. \cite{Arenc}),
which we intend to be defined as the energy function which is preserved whenever $\mu=0$.
Let us write the normal form equations \equ{NFf} using the following compact notation as in \equ{NFcompact0}:
\beqa{NFcompact3}
\dot X &=& \omega(Y)+\varepsilon p_Y^{(\leq K)}(Y,X,t;\varepsilon)+\mu s^{(\leq K)}(Y,X,t;\varepsilon,\mu)\nonumber\\
&+&F^{(>K)}(Y,X,t;\varepsilon,\mu)+F_{N+1}(Y,X,t;\varepsilon,\mu)\nonumber\\
\dot Y &=& -\varepsilon p_X^{(\leq K)}(Y,X,t;\varepsilon)+G^{(>K)}(Y,X,t;\varepsilon,\mu)+G_{N+1}(Y,X,t;\varepsilon,\mu)\nonumber\\
\dot U &=& -\varepsilon p_t^{(\leq K)}(Y,X,t;\varepsilon)+H^{(>K)}(Y,X,t;\varepsilon,\mu)+H_{N+1}(Y,X,t;\varepsilon,\mu)\ ,
\eeqa
where we have indicated also the dependence on the parameters and
we have identified the functions as follows:
\beqano
\varepsilon p_Y^{(\leq K)}(Y,X,t;\varepsilon)&\equiv&\sum_{i=1}^N\Omega_{i0}(Y)\varepsilon^i+
\sum_{i=1}^N F_{i0}^{(3,r,\leq K,\leq N)}(Y,X,t)\varepsilon^i\nonumber\\
\varepsilon p_X^{(\leq K)}(Y,X,t;\varepsilon)&\equiv&-\sum_{i=1}^N G_{i0}^{(3,r,\leq K,\leq N)}(Y,X,t)\varepsilon^i\nonumber\\
\varepsilon p_t^{(\leq K)}(Y,X,t;\varepsilon)&\equiv&-\sum_{i=1}^N H_{i0}^{(3,r,\leq K,\leq N)}(Y,X,t)\varepsilon^i\nonumber\\
\mu s^{(\leq K)}(Y,X,t;\varepsilon,\mu)&\equiv& \sum_{j=1}^N \sum_{i=0}^{N-j} F_{ij}^{(3,r,\leq K,\leq N)}(Y,X,t)\varepsilon^i\mu^j+
\sum_{j=1}^N \sum_{i=0}^{N-j} \Omega_{ij}(Y)\varepsilon^i\mu^j\ .
\eeqano
For $\mu=0$ equations \equ{NFcompact3} reduce to
\beqa{NFcompactmu}
\dot X &=& \omega(Y)+\varepsilon p_Y^{(\leq K)}(Y,X,t;\varepsilon)+
F^{(>K)}(Y,X,t;\varepsilon,0)+F_{N+1}(Y,X,t;\varepsilon,0)\nonumber\\
\dot Y &=& -\varepsilon p_X^{(\leq K)}(Y,X,t;\varepsilon)+G^{(>K)}(Y,X,t;\varepsilon,0)+G_{N+1}(Y,X,t;\varepsilon,0)\nonumber\\
\dot U &=& -\varepsilon p_t^{(\leq K)}(Y,X,t;\varepsilon)+H^{(>K)}(Y,X,t;\varepsilon,0)+H_{N+1}(Y,X,t;\varepsilon,0)\ .
\eeqa
Due to the Hamiltonian character of the equations of motion for $\mu=0$, there exist vector functions $A^{(>K)}$, $B_{N+1}$, such that
\beqano
A_Y^{(>K)}=F^{(>K)}(Y,X,t;\varepsilon,0)\ ,\quad A_X^{(>K)}=-G^{(>K)}(Y,X,t;\varepsilon,0)\ ,\quad A_t^{(>K)}=-H^{(>K)}(Y,X,t;\varepsilon,0)\nonumber\\
B_{N+1,Y}=F_{N+1}(Y,X,t;\varepsilon,0)\ ,\quad B_{N+1,X}=-G_{N+1}(Y,X,t;\varepsilon,0)\ ,\quad B_{N+1,t}=-H_{N+1}(Y,X,t;\varepsilon,0)\ ,
\eeqano
so that we can recognize \equ{NFcompactmu} as
Hamilton's equations associated to the following Hamiltonian function in the extended phase space with $\dot t=1$:
$$
{\cal H}(Y,X,U,t)=h_{00}(Y)+U+\varepsilon p^{(\leq K)}(Y,X,t)+A^{(>K)}(Y,X,t)+B_{N+1}(Y,X,t)\ ,
$$
where $h_{00}$ is such that
${{\partial h_{00}(Y)}\over {\partial Y}}=\omega(Y)$. Let us fix the energy level
${\cal H}=E$ for some real constant $E$; taking into account the complete equations \equ{NFcompact3}, we obtain that the
variation of $E$ for $\mu\not=0$ is given by (for simplicity we omit the arguments):
\beq{53bis3}
{dE\over dt}=\varepsilon\mu\, p_X^{(\leq K)}s^{(\leq K)}+C_{N+1}+D^{(>K)}\ ,
\eeq
with
\beqano
C_{N+1}&\equiv&\omega(Y)G_{N+1}+\varepsilon p_Y^{(\leq K)} G_{N+1}+\varepsilon p_X^{(\leq K)} F_{N+1}\nonumber\\
&+&(A_Y^{(>K)}+B_{N+1,Y})G_{N+1}+B_{N+1,Y}(-\varepsilon p_X^{(\leq K)}+G^{(>K)})\nonumber\\
&+&(A_X^{(>K)}+B_{N+1,X})F_{N+1}+B_{N+1,X}(\omega+\varepsilon p_Y^{(\leq K)}+\mu s^{(\leq K)}+F^{(>K)})
+H_{N+1}+B_{N+1,t}\nonumber\\
D^{(>K)}&\equiv&\omega(Y)G^{(>K)}+\varepsilon p_Y^{(\leq K)} G^{(>K)}+\varepsilon p_X^{(\leq K)} F^{(>K)}+H^{(>K)}\nonumber\\
&-&\varepsilon p_X^{(\leq K)} A_Y^{(>K)}+A_Y^{(>K)} G^{(>K)}+\omega A_X^{(>K)}+\varepsilon p_Y^{(\leq K)} A_X^{(>K)}\nonumber\\
&+&\mu s^{(\leq K)}\,A_X^{(>K)}+A_X^{(>K)} F^{(>K)}+A_t^{(>K)}\ ,
\eeqano
where now the functions $F_{N+1}$, $G_{N+1}$, $H_{N+1}$, $F^{(>K)}$, $G^{(>K)}$, $H^{(>K)}$ depend on $(Y,X,t;\varepsilon,\mu)$.
Denoting by $\Delta E\equiv E(t)-E(0)$, we obtain
$$
|\Delta E|\geq |\Delta h_{00}+\Delta U|-\Big(\varepsilon\|\Delta p^{(\leq K)}\|_{R_0,S_0}+
\|\Delta A^{(>K)}\|_{R_0,S_0}+\|\Delta B_{N+1}\|_{R_0,S_0}\Big)\ ,
$$
where $\Delta h_{00}+\Delta U\equiv h_{00}(Y(t))-h_{00}(Y(0))+U(t)-U(0)$ and similarly for the other quantities.
Recalling \equ{53bis3} and setting $h_0\equiv h_{00}+U$, we get
$$
|\Delta h_0|\leq |\Delta E|+\varepsilon\|\Delta p^{(\leq K)}\|_{R_0,S_0}+
\|\Delta A^{(>K)}\|_{R_0,S_0}+\|\Delta B_{N+1}\|_{R_0,S_0}\ ,
$$
where
$$
|\Delta E|\leq|{{dE}\over {dt}}|\,t \leq \Big(\varepsilon\mu \|p_X^{(\leq K)}\|_{R_0,S_0}
\|s^{(\leq K)}\|_{R_0,S_0}+\|C_{N+1}\|_{R_0,S_0}+\|D^{(>K)}\|_{R_0,S_0}\Big)\,t\ .
$$
We denote by $m$ an upper bound on the Hessian of $h_{00}(y)$ and let $\tilde m$ be an upper bound of the Hessian in
the normalized variables, which we can define as $\tilde m\equiv m+\|{{\partial^3 h_{00}}\over {\partial y^3}}\|_{r_0}
\,\|D^{(y,N)}\|_{r_0,s_0}$, having expressed the link between new and old variables as $Y=y+D^{(y,N)}(y,x,t)$.
Then, we have:
$$
\sup_{Y\in C_{r_0}(A)} \|{{\partial^2 h_{00}(Y)}\over {\partial Y^2}}\|\leq \tilde m\ .
$$
Assume that the frequency $\omega_e(y)\equiv(\omega(y),1)$ is close to exact $\Lambda$--resonances
(compare with \cite{poe1}) by a quantity $\delta>0$, namely if
$R_\Lambda\equiv \{\Omega\in{\R}^{\ell+1}:\ \Omega\cdot n=0$ for all $n\in\Lambda\}$, then
$\min_{\Omega\in R_\Lambda}\|\omega_e(y)-\Omega\|\leq \delta$.
Setting $Z\equiv(Y,U)$, assume that
$|Z(t)-Z(0)|\leq r$ for some $r>0$ with $\delta+\tilde mr\leq R_0$.
Let $\Pi_\Lambda$ be the orthogonal projection on $\Lambda$; by the mean value theorem we obtain
\beqano
|\omega_e\cdot \Pi_\Lambda\Delta Z|&\leq& \|\Pi_\Lambda\omega_e\|_{R_0}\,\|P\Delta Z\|\leq (\delta+\tilde mr)\,\|\Delta Z\|\nonumber\\
|\omega_e\cdot (Id.-\Pi_\Lambda)\Delta Z|&\leq& t\ \|\omega_e\|_{R_0}\ \left(\|G^{(>K)}\|_{R_0,S_0}+\|G_{N+1}\|_{R_0,S_0}+
\|H^{(>K)}\|_{R_0,S_0}+\|H_{N+1}\|_{R_0,S_0}\right)\ .
\eeqano
Moreover:
$$
\Delta h_0=\omega_e\cdot \Delta Z+\int_0^1(1-s){{\partial^2h_{0}(Z(s))}\over{\partial Z^2}}\Delta Z\cdot\Delta Z\ ds\ ,
$$
so that, in the region where the $M$--convexity \equ{QC} holds, one has
$$
{\tilde M\over 2}\|\Delta Z\|^2\leq |\omega_e\cdot \Delta Z|+|\Delta h_0|\ ,
$$
where (similarly to $\tilde  m$) we can set $\tilde M\equiv M-\|{{\partial^3 h_{00}}\over {\partial y^3}}\|_{r_0}
\,\|D^{(y,N)}(y,x,t)\|_{r_0,s_0}$, so that one has
$$
{{\partial^2h_{0}(Z)}\over{\partial Z^2}}v\cdot v\geq \tilde M\|v\|^2\ ,\qquad \forall v\in{\R}^{\ell+1}\ .
$$
Finally, we have
\beqano
{\tilde M\over 2}\|\Delta Z\|^2&\leq& (\delta+\tilde mr)\|\Delta Z\|+t\|\omega_e\|_{R_0}
\Big(\|G^{(>K)}\|_{R_0,S_0}+\|G_{N+1}\|_{R_0,S_0}\nonumber\\
&+&\|H^{(>K)}\|_{R_0,S_0}+\|H_{N+1}\|_{R_0,S_0}\Big)+\varepsilon\mu \|p_X^{(\leq K)}\|_{R_0,S_0}\,
\|s^{(\leq K)}\|_{R_0,S_0}\,t+\|C_{N+1}\|_{R_0,S_0}\,t\nonumber\\
&+&\|D^{(>K)}\|_{R_0,S_0}\,t+\varepsilon\|\Delta p^{(\leq K)}\|_{R_0,S_0}+
\|\Delta A^{(>K)}\|_{R_0,S_0}+\|\Delta B_{N+1}\|_{R_0,S_0}\  ,
\eeqano
which gives a bound on the norm of $\Delta Z$. Notice that $\|G^{(>K)}\|$, $\|H^{(>K)}\|$ and $\|\Delta A^{(>K)}\|$ are
of order of $e^{-K\tau_0}$, namely of order $\lambda^N$ once we set $N$ such that $N=[K\tau_0/|\log\lambda|]$. We finally define the
constants $C_1$, $C_2$, $C_3$, $C_4$ such that
\beqa{CC}
\|\omega_e\|_{R_0} \Big(\|G^{(>K)}\|_{R_0,S_0}&+&\|G_{N+1}\|_{R_0,S_0}+
\|H^{(>K)}\|_{R_0,S_0}+\|H_{N+1}\|_{R_0,S_0}\Big)\nonumber\\
+\|C_{N+1}\|_{R_0,S_0}+\|D^{(>K)}\|_{R_0,S_0}&\leq& C_1\lambda^N\nonumber\\
\|p_X^{(\leq K)}\|_{R_0,S_0}\|s^{(\leq K)}\|_{R_0,S_0}&\leq& C_2\nonumber\\
\|\Delta p^{(\leq K)}\|_{R_0,S_0}&\leq& C_3\nonumber\\
\|\Delta A^{(>K)}\|_{R_0,S_0}+\|\Delta B_{N+1}\|_{R_0,S_0}&\leq& C_4\lambda^N\ .
\eeqa
With this setting we obtain:
$$
{\tilde M\over 2}\|\Delta Z\|^2\leq (\delta+\tilde mr)\|\Delta Z\|+C_1\,\lambda^N\,t+C_2\,\varepsilon\mu\,t+C_3\,\varepsilon +C_4\,\lambda^N\ .
$$
Based on the above formula and on $\|\Delta Y\|\leq \|\Delta Z\|$, we can draw the following conclusions:
\begin{enumerate}
    \item for some $\rho_1>0$, $\|\Delta Y\|\leq\rho_1$ for $t$ of the order of $\lambda^{-N}$ if\\
$C_2=0$, namely if $\|p_X^{(\leq K)}\|_{R_0,S_0} \|s^{(\leq K)}\|_{R_0,S_0}=0$, i.e.
    either $\|p_X^{(\leq K)}\|_{R_0,S_0}=0$ or $\|s^{(\leq K)}\|_{R_0,S_0}=0$;
    \item for some $\rho_2>0$, $\|\Delta Y\|\leq\rho_2$ for $t$ of the order of the minimum
between $\lambda^{-N}$ and
    $(\varepsilon\mu)^{-1}$ if $C_2\not=0$, i.e. $\|p_x^{(\leq K)}\|_{R_0,S_0} \|s^{(\leq K)}
\|_{R_0,S_0}\not=0$.
   \end{enumerate}

\noindent
The two cases correspond to items $i)$, $ii)$ of the statement of the Theorem. More precisely, let us start with
the case $C_2=0$, i.e. $p_X^{(\leq K)}=0$ or $s^{(\leq K)}=0$.
Assuming that $\delta+\tilde mr  <\alpha \tilde Mr$,
$C_3\varepsilon+C_4 \lambda^N< \beta \tilde Mr^2$, $t< {1\over {C_1\lambda^N}}\gamma \tilde Mr^2$ with $C_1>0$,
for some positive constants $\alpha$, $\beta$, $\gamma$,
under the assumption that $\|\Delta Z\|\leq r$ for some $r>0$, we obtain
$$
{\tilde M\over 2}\|\Delta Z\|^2< (\alpha+\beta+\gamma) \tilde Mr^2\ ,
$$
namely
\beq{rho1}
\|\Delta Y\|\leq\|\Delta Z\|< \sqrt{2(\alpha+\beta+\gamma)}\, r\equiv \rho_1
\eeq
with $\rho_1\leq R_0$.
Taking into account \equ{yY} and \equ{Pi}, one obtains item $i)$
of the Theorem, namely
\beq{T1}
\|y(t)-y(0)\|\leq 2C_p\,\lambda+\rho_1\ \quad {\rm for}\ t\leq T_1\equiv C_0\, e^{K\tau_0}
\eeq
with $C_0\equiv(\gamma \tilde Mr^2)/C_1$.

\noindent
Concerning item $ii)$, since $C_2\not=0$ let $\sigma>0$ be such that for
$t< \min({1\over {C_1\lambda^N}}\gamma \tilde Mr^2,{1\over {C_2\varepsilon\mu}}\,\sigma \tilde Mr^2)$,
one has
$$
{\tilde M\over 2}\|\Delta Z\|^2<(\alpha+\beta+\gamma+\sigma) \tilde Mr^2\ ,
$$
namely
\beq{rho2}
\|\Delta Y\|\leq \|\Delta Z\|< \sqrt{2(\alpha+\beta+\gamma+\sigma)}\, r\equiv\rho_2\ ,
\eeq
with $\rho_2\leq R_0$. According to \equ{yY} and \equ{Pi}, we obtain that
\beq{T2}
\|y(t)-y(0)\|\leq 2C_p\,\lambda+\rho_2\ \quad {\rm for}\ t\leq T_2\equiv \min({{C_0}\over {\lambda^N}},\, {C_0'\over {\varepsilon\mu}})
\eeq
with $C_0\equiv(\gamma \tilde Mr^2)/C_1$, $C_0'\equiv (\sigma \tilde Mr^2)/C_2$.

\noindent
When the $M$--convexity condition is violated (i.e. the second condition in \equ{QC}), by the assumption of quasi--convexity
the first inequality in \equ{QC} must hold. Let $T$ be either $T_1$ or $T_2$ as in \equ{T1}, \equ{T2} with
$\|\omega_e\|_{R_0}$ replaced by
$$
C_\omega\equiv \sup_{\|Y-Y_0\|\leq \rho}\|\omega_e(Y)\|
$$
with $\rho$ being $\rho_1$ or $\rho_2$ as in \equ{rho1}, \equ{rho2}. Repeating the same argument as in \cite{poe1}, we assume that there exists an escape
time $T_e$ such that $\|\Delta Z\|=\rho$ and we show that this implies the inequality
$$
|\omega_e(Y(s))\cdot \Delta Z|\leq L\,\|\Delta Z\|\qquad \forall\ 0\leq s\leq 1\ .
$$
Then, using the same argument as for the convex case, we conclude that $\|\Delta Z\|<\rho$, thus providing a contradiction.
As before we have:
\beqano
|\omega_e(Y(s))\cdot P\Delta Z|&\leq& (\delta+\tilde mr)\,\|\Delta Z\|\nonumber\\
|\omega_e(Y(s))\cdot (Id.-P)\Delta Z|&\leq& C_1\lambda^N\, T\ .
\eeqano
Then, we have
$$
|\omega_e(Y(s))\cdot \Delta Z|\leq (\delta+\tilde mr)\,\|\Delta Z\|+C_1\lambda^N\, T\ .
$$
If $T=T_1$, $\rho=\rho_1$, we obtain
\beqano
|\omega_e(Y(s))\cdot \Delta Z|&\leq&(\delta+\tilde mr)\,\|\Delta Z\|+C_1\lambda^N\, T_1\nonumber\\
&<&\alpha \tilde M r\rho_1+\gamma\tilde M r^2\nonumber\\
&=&({\alpha\over {\sqrt{2(\alpha+\beta+\gamma)}}}+{\gamma\over {2(\alpha+\beta+\gamma)}})\tilde M \rho_1^2\nonumber\\
&\leq&\ L\|\Delta Z\|\ ,
\eeqano
if
$$
\rho_1\leq {L\over {({\alpha\over {\sqrt{2(\alpha+\beta+\gamma)}}}+{\gamma\over {2(\alpha+\beta+\gamma)}})\tilde M}}\ .
$$
If $T=T_2$, $\rho=\rho_2$, assume that $T=(\sigma \tilde M r^2)/(C_2\varepsilon\mu)$ (otherwise we recover the case $T=T_1$). Then,
$$
|\omega_e(Y(s))\cdot \Delta Z|\leq (\delta+\tilde mr)\,\|\Delta Z\|+C_1\lambda^N\, {{\sigma \tilde M r^2}\over {C_2\varepsilon\mu}}\ .
$$
Being
$$
{{C_1\lambda^N}\over {C_2\varepsilon\mu}}\ \leq\ {\gamma\over\sigma}\ ,
$$
we obtain:
\beqano
|\omega_e(Y(s))\cdot \Delta Z|&<& (\delta+\tilde mr)\,\|\Delta Z\|+{\gamma\over\sigma}\ \sigma\tilde M r^2\nonumber\\
&<&\alpha \tilde M r\|\Delta Z\|+\gamma\tilde M r^2\nonumber\\
&=&({\alpha\over {\sqrt{2(\alpha+\beta+\gamma+\sigma)}}}+{\gamma\over {2(\alpha+\beta+\gamma+\sigma)}})\tilde M \rho_2^2\nonumber\\
&\leq& L\|\Delta Z\|\ ,
\eeqano
which is satisfied if the following condition holds:
$$
\rho_2\leq {L\over {({\alpha\over {\sqrt{2(\alpha+\beta+\gamma+\sigma)}}}+{\gamma\over {2(\alpha+\beta+\gamma+\sigma)}})\tilde M}}\ .
$$
$\Box$

\vskip.1in

\noindent
\bf Remark 5. \rm Since we do not claim the result for any $y_0\in A$, but only locally under the conditions \equ{res} and \equ{omega0},
we do not need to cover the whole phase space and therefore we do not need the analysis of the geography of the resonances,
as it is usually done (see, e.g., \cite{poe1}).

\section{Applications of the normal forms}\label{sec:examples}
As we have seen in the Theorem, the stability time depends on the expressions of the terms $p_X^{(\leq K)}$ and $s^{(\leq K)}$ appearing in
the normal form equations, which represent, respectively, the conservative resonant part of the action variables and
the dissipative resonant part pertaining to the angles, including the contribution of the modified frequency.
In this Section we analyze several different examples,
which well represent all possible situations which can be obtained with different choices of $p_X^{(\leq K)}$ and $s^{(\leq K)}$.
We illustrate these models with a twofold goal: to provide examples of cases $i)$ and $ii)$ of the Theorem
and to illustrate an explicit evaluation of the resonant normal form.
Since we do not aim to obtain stability estimates,
we limit ourselves to the computation of the normal form in the non--extended phase space, i.e. in the variables $x$ and $y$ only.
The experiments performed in this Section will be validated by the theoretical results of Section~\ref{sec:APP},
where the estimates of the Theorem will be applied, showing linear as well as exponential stability times.

\vskip.1in

\noindent
All examples considered in the forthcoming Sections~\ref{EFLD}--\ref{EPXS2} will have the following simple form:
\beqano
\dot x&=&y+\mu\,f_{01}(x,t)\nonumber\\
\dot y&=&-\varepsilon h_{10,x}(x,t)-\mu(y-\eta)\ ,
\eeqano
where $f_{01}$ and $h_{10}$ are periodic functions.
In this case it is easy to decide which of the conditions $i)$ or $ii)$ of the Theorem are satisfied.
Since we shall not need to consider Fourier modes less or equal, or greater than $K$, we drop the superscript
by writing $p_X$, $p_Y$, $s$ in place of $p_X^{(\leq K)}$, $p_Y^{(\leq K)}$, $s^{(\leq K)}$. Then, we can state
that $p_X=0$ whenever the resonant part of $h_{10,x}$ is zero, otherwise $p_X$ is different from zero.
Concerning the function $s$, we can state that if the resonant part of $f_{01}$ is not zero as well
as if products of the form $(f_{01})^m(h_{10,x})^n$ with $0<m+n\leq N$ generate resonant terms
or zero average terms of order $\mu$, then
the function $s$ is different from zero. If the products $(f_{01})^m(h_{10,x})^n$ with $0<m+n\leq N$ do not
generate resonant terms or zero average terms of order $\mu$, then $s=0$ up to the order $N$.

\subsection{Linear stability: case $p_X\not=0$, $s\not=0$}\label{EFLD}

\begin{figure}[htp]
\includegraphics[width=4.8cm]{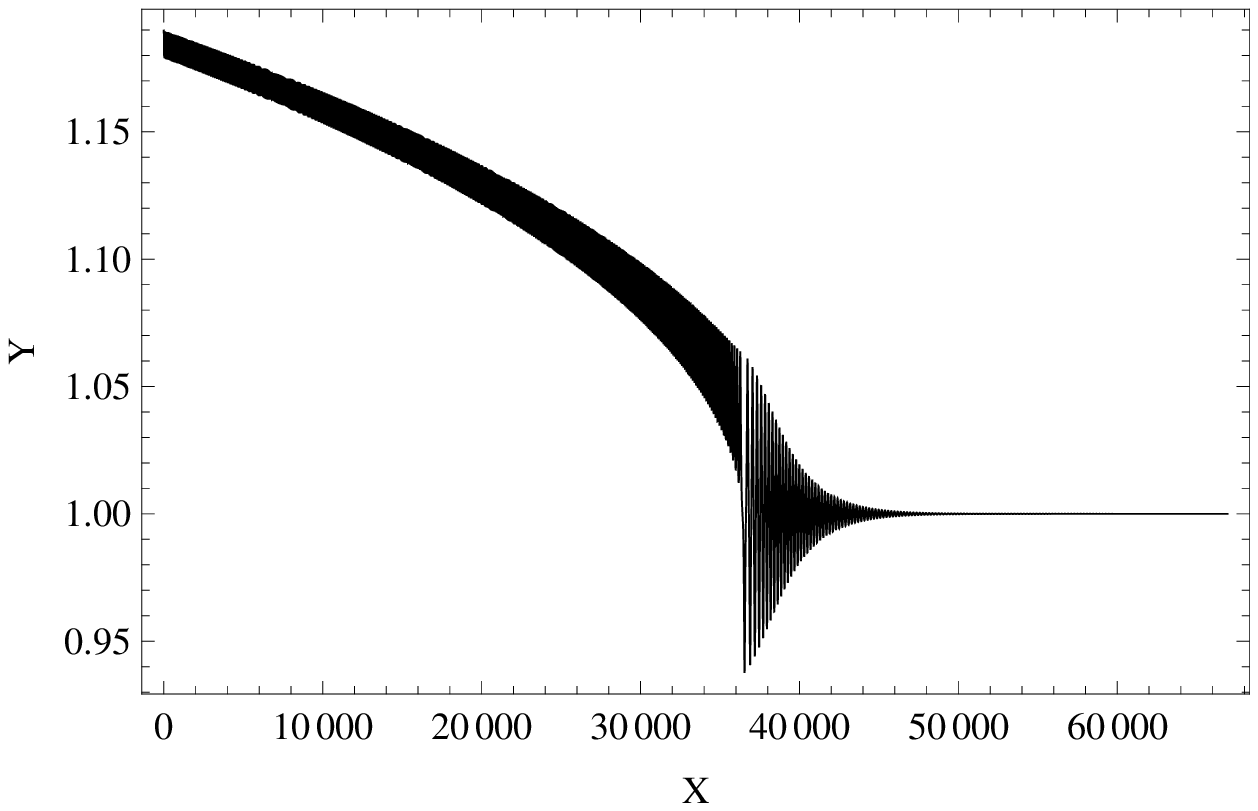}
\hglue-0.01cm
\includegraphics[width=4.8cm]{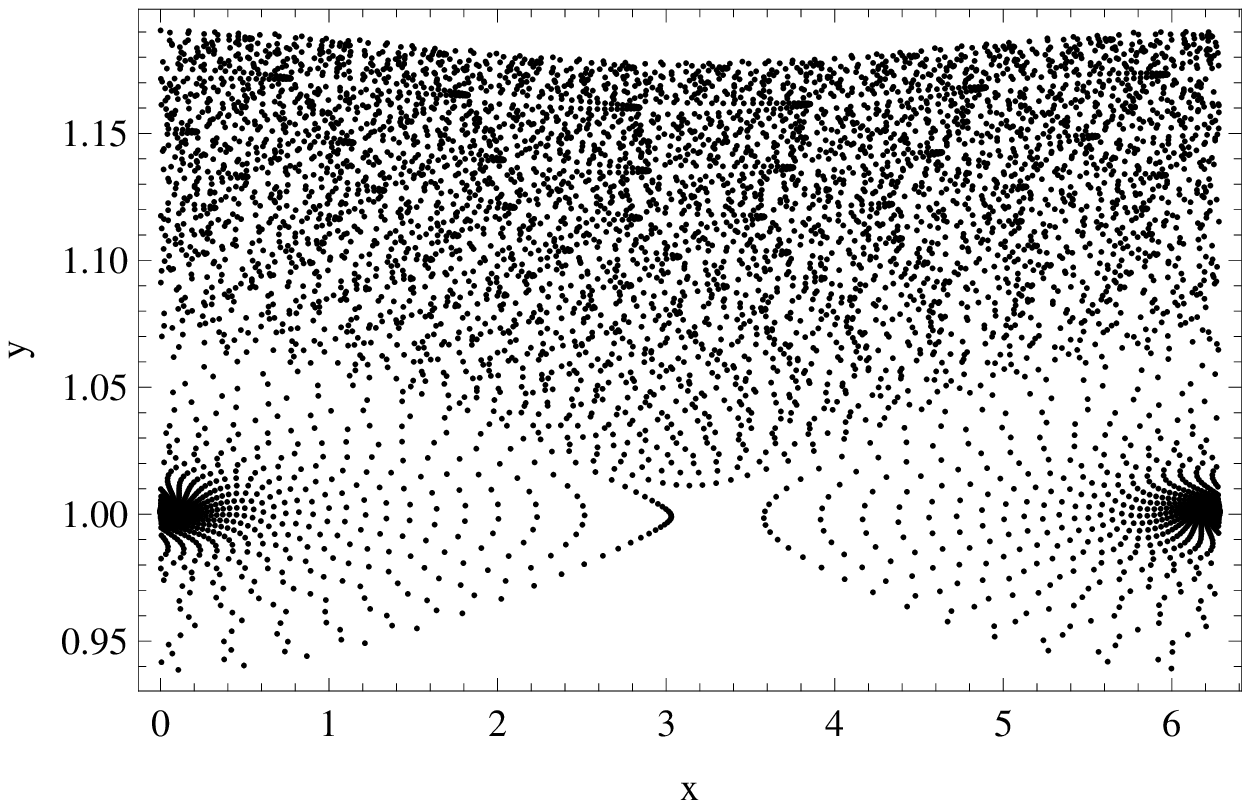}
\hglue-0.01cm
\includegraphics[width=5.3cm]{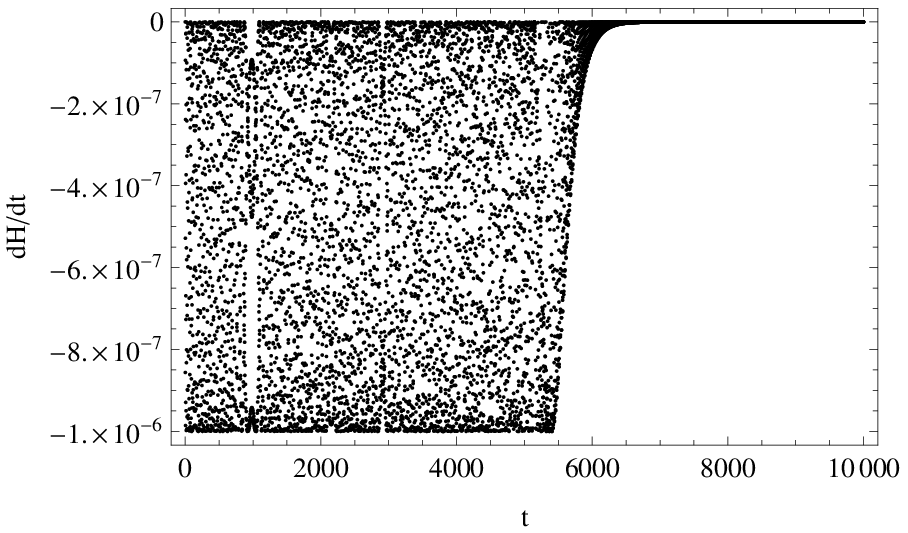}
\caption{Case $p_X\not=0$, $s\not=0$ associated to \equ{e19}
for $\varepsilon=10^{-3}$, $\mu=10^{-3}$ and for the initial conditions
$X(0)=0$, $Y(0)=1+6\sqrt{\varepsilon}$. Left: the lift of the normal form variables
$(X,Y)$ to the universal coverage. Middle: the trajectory in the original variables $(x,y)$.
Right: the variation of the derivative of the normalized Hamiltonian.}\label{exa1-orb}
\end{figure}

We consider the one--dimensional, time--dependent vector field given by
\beqa{e19}
\dot x&=&y-\mu (\sin (x-t)+\sin (x)) \nonumber \\
\dot y&=&-\varepsilon (\sin (x-t)+\sin (x)) -\mu  (y-\eta) \ .
\eeqa
Following the calculations of the proof of the Resonant Normal Form Lemma, the conservative transformation
up to second order is defined by
\beqano
\psi _{10}(\tilde y,\tilde x,t)&=&\frac{\sin (\tilde x)}{\tilde y} \nonumber \\
\psi _{20}(\tilde y,\tilde x,t)&=&\frac{\sin (2 \tilde x-t)}{2 \tilde y^2 (2 \tilde y-1)}-\frac{\sin (t)}{2 \tilde y^2}-\frac{\sin (2 \tilde x)}{8 \tilde y^3} \ ,
\eeqano
while the dissipative contribution is given by
\beqano
\beta _{01}(Y,X,t)&=&0\nonumber \\
\alpha _{01}(Y,X,t)&=&-\frac{\cos (X)}{Y}\nonumber \\
\beta _{11}(Y,X,t)&=&-\frac{\sin (2 X)}{4 Y^2}+\frac{\sin (t)}{Y} \nonumber \\
\alpha _{11}(Y,X,t)&=&-\frac{\cos (2 X-t)}{2 Y^2 (2 Y-1)}-\frac{(2 Y+1) \cos (t)}{2 Y^2}+\frac{\cos (2 X)}{8 Y^3}  \nonumber \\
\beta _{02}(Y,X,t)&=&0 \nonumber \\
\alpha _{02}(Y,X,t)&=&\frac{\sin (t)}{Y}+\frac{\sin (2 X)}{4 Y^2}\ .
\eeqano
By choosing
$$
\eta(Y) =Y+\frac{\varepsilon }{2 Y} + O_3(\varepsilon,\mu)  \ ,
$$
the normal form equations become
\beqano
\dot X&=&Y-\frac{\varepsilon ^2}{2 Y^3}-\frac{\mu ^2}{2 Y}-\mu  \sin (X-t) + O_3(\varepsilon,\mu) \nonumber \\
\dot Y&=&-\varepsilon  \sin (X-t) + O_3(\varepsilon,\mu) \ ,
\eeqano
where we recognize that $p_X(Y,X,t)=\sin(X-t)$, $p_Y(Y,X,t)=-\varepsilon/(2Y^3)$, $s(Y,X,t)=-\sin(X-t)-\frac{\mu}{2 Y}$.
The Hamiltonian function in the extended phase space with $U$ conjugated to time, associated to the normalized equations
for $\mu=0$, is given by
$${\cal H}(Y,X,U,t)=\frac{Y^2}{2}+U+\frac{\varepsilon ^2}{4 Y^2}-\varepsilon  \cos (X-t)+O_3(\varepsilon,\mu) \ .$$
Replacing the normalized equations into the total derivative of ${\cal H}$, one gets
\beq{der}
\frac{d{{\cal H}(Y,X,U,t)}}{dt}=-\frac{1}{2} \mu  \varepsilon \Big(1-\cos (2 X-2t)\Big)+O_3(\varepsilon,\mu) \ .
\eeq
A typical orbit is shown in Figure~\ref{exa1-orb}, where we integrate the normal form equations for $\varepsilon=10^{-3}$
and $\mu=10^{-3}$ with initial conditions $X(0)=0$ and $Y(0)=1+6\sqrt{\varepsilon}$. The left panel of Figure \ref{exa1-orb} shows the
lift of $(X,Y)$ to the universal coverage, while the middle panel shows the orbit back--transformed to the old variables $(x,y)$.
The dynamics starts on a rotational regime and drifts downwards; then it spirals along librational invariant curves until
reaching the attractor. The right panel of Figure~\ref{exa1-orb} provides the variation of the derivative of the
normal form Hamiltonian, which tends to zero as the orbit reaches the attractor.
The behavior is justified by \equ{der} as the resonance is approached.

\subsection{Linear stability at higher orders: case $p_X\not=0$, $s\not=0$}\label{ESLD}

We consider the vector field
\beqa{e20}
\dot x&=&y-\mu \sin (x) \nonumber \\
\dot y&=&-\varepsilon (\sin (x-t)+\sin (x)) -\mu  (y-\eta)\ .
\eeqa
The conservative normal form is defined by
\beqano{}
\psi _{10}(\tilde y,\tilde x,t)&=&\frac{\sin (\tilde x)}{\tilde y} \nonumber \\
\psi _{20}(\tilde y,\tilde x,t)&=&\frac{\sin (2 \tilde x-t)}{2 \tilde y^2 (2 \tilde y-1)}-\frac{\sin (t)}{2 \tilde y^2}-\frac{\sin (2 \tilde x)}{8 \tilde y^3}  \ .
\eeqano
The dissipative transformation becomes:
\beqano{}
\beta _{01}(Y,X,t)&=&0 \nonumber \\
\alpha _{01}(Y,X,t)&=&-\frac{\cos (X)}{Y} \nonumber \\
\beta _{11}(Y,X,t)&=&\frac{\sin (2 X-t)}{2 Y (2 Y-1)}+\frac{\sin (t)}{2 Y}-\frac{\sin (2 X)}{4 Y^2} \nonumber \\
\alpha _{11}(Y,X,t)&=&\frac{(1-3 Y) \cos (2 X-t)}{2 Y^2 (2 Y-1)^2}+\frac{(1-Y) \cos (t)}{2 Y^2}+\frac{\cos (2 X)}{8 Y^3} \nonumber \\
\beta _{02}(Y,X,t)&=&0 \nonumber \\
\alpha _{02}(Y,X,t)&=&\frac{\sin (2 X)}{4 Y^2} \ .
\eeqano
Higher normal form terms associated to \equ{e20} can be obtained in a similar way.
For this model resonant terms occur at higher orders; for this reason we provide
the following third order normal form equations:
\beqano{}
\dot X&=&Y-\frac{\varepsilon ^2}{2 Y^3}-\frac{\mu ^2}{2 Y} + \frac{\varepsilon ^3 (2-5 Y) \cos (X-t)}{2 (1-2 Y)^2 Y^5}
-\frac{\varepsilon ^2 \mu  \left(6 Y^2+2 Y-1\right) \sin (X-t)}{2 (1-2 Y)^2 Y^3}\nonumber \\
&+&\frac{1}{2} \varepsilon  \mu ^2 \left(\frac{1}{Y^2}-\frac{2 Y}{(1-2 Y)^2}\right) \cos
   (X-t) + O_4(\varepsilon,\mu) \nonumber \\
\dot Y&=&-\varepsilon  \sin (X-t) +\frac{\varepsilon ^3 \sin (X-t)}{8 Y^5-4 Y^4} + O_4(\varepsilon,\mu) \ .
\eeqano
The  $\eta$ is the same as in the previous Section, the Hamiltonian associated to the normal form in
the extended phase space is given by
\beqno
{\cal H}(Y,X,U,t)=\frac{Y^2}{2}+U+\frac{\varepsilon ^2}{4 Y^2}-\varepsilon  \cos (X-t)+\frac{\varepsilon ^3 \cos (X-t)}{4 Y^4 (2 Y-1)} +O_4(\varepsilon,\mu) \ ,
\eeqno
while the derivative of the Hamiltonian becomes
\beqno
\frac{d{\cal H}(Y,X,U,t)}{dt}=-\frac{\varepsilon  \mu ^2 \sin (X-t)}{2 Y} +O_4(\varepsilon,\mu)\ .
\eeqno
The normal form produces a resonant term at third order. As a consequence, we observe a drift of the action variables,
but on longer time scales.

\begin{figure}[htp]
\includegraphics[width=4.8cm]{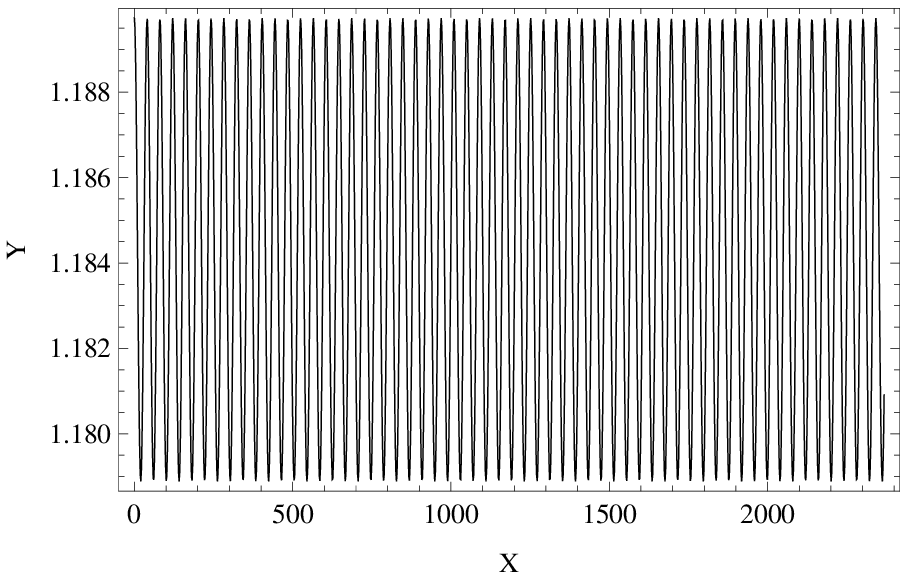}
\hglue-0.01cm
\includegraphics[width=4.8cm]{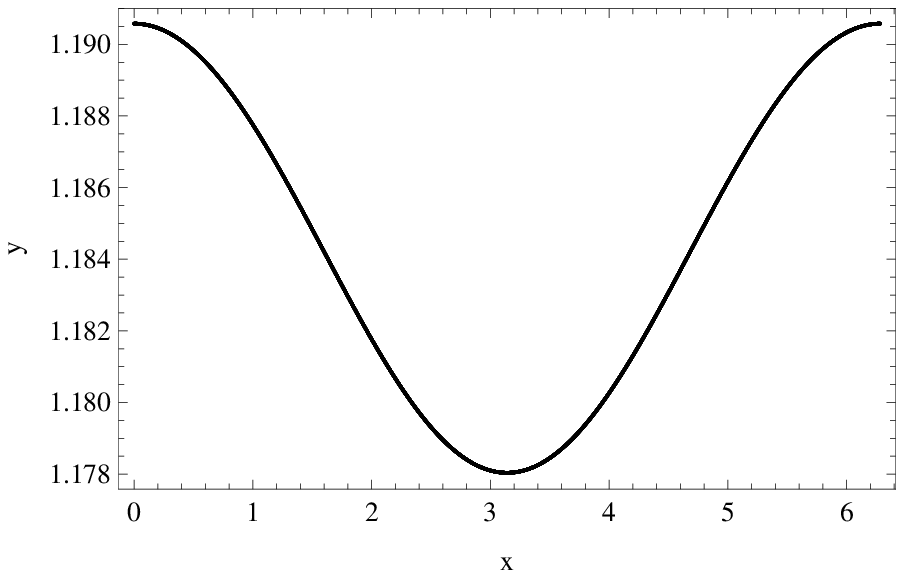}
\hglue-0.01cm
\includegraphics[width=5.3cm]{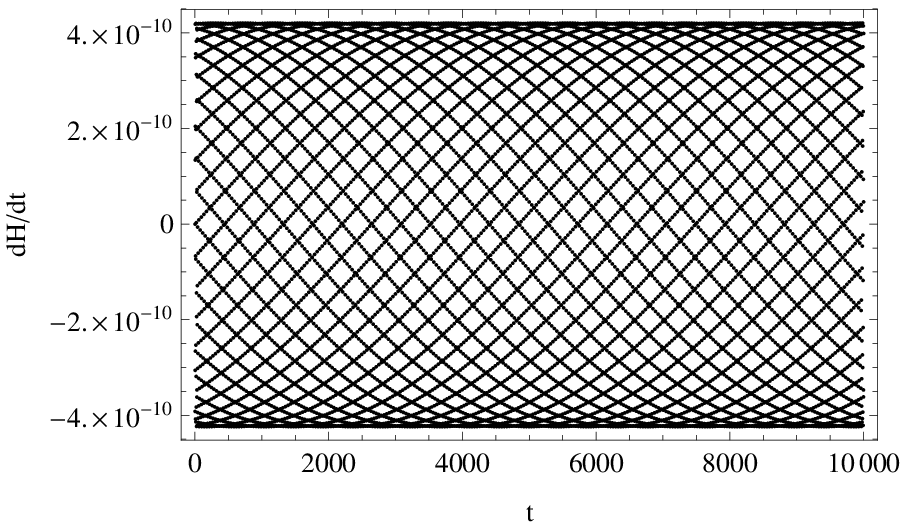}
\caption{Case $p_X\not=0$, $s\not=0$ associated to \equ{e20}
for $\varepsilon=10^{-3}$, $\mu=10^{-3}$ and for the initial conditions
$X(0)=0$, $Y(0)=1+6\sqrt{\varepsilon}$. Left: the lift of the normal form variables
$(X,Y)$ to the universal coverage. Middle: the trajectory in the original variables $(x,y)$.
Right: the variation of the derivative of the normalized Hamiltonian.}\label{exa2-orb}
\end{figure}

\subsection{Exponential stability: case $p_X=0$, $s\neq0$}\label{EPXS1}
\begin{figure}[htp]
\includegraphics[width=4.8cm]{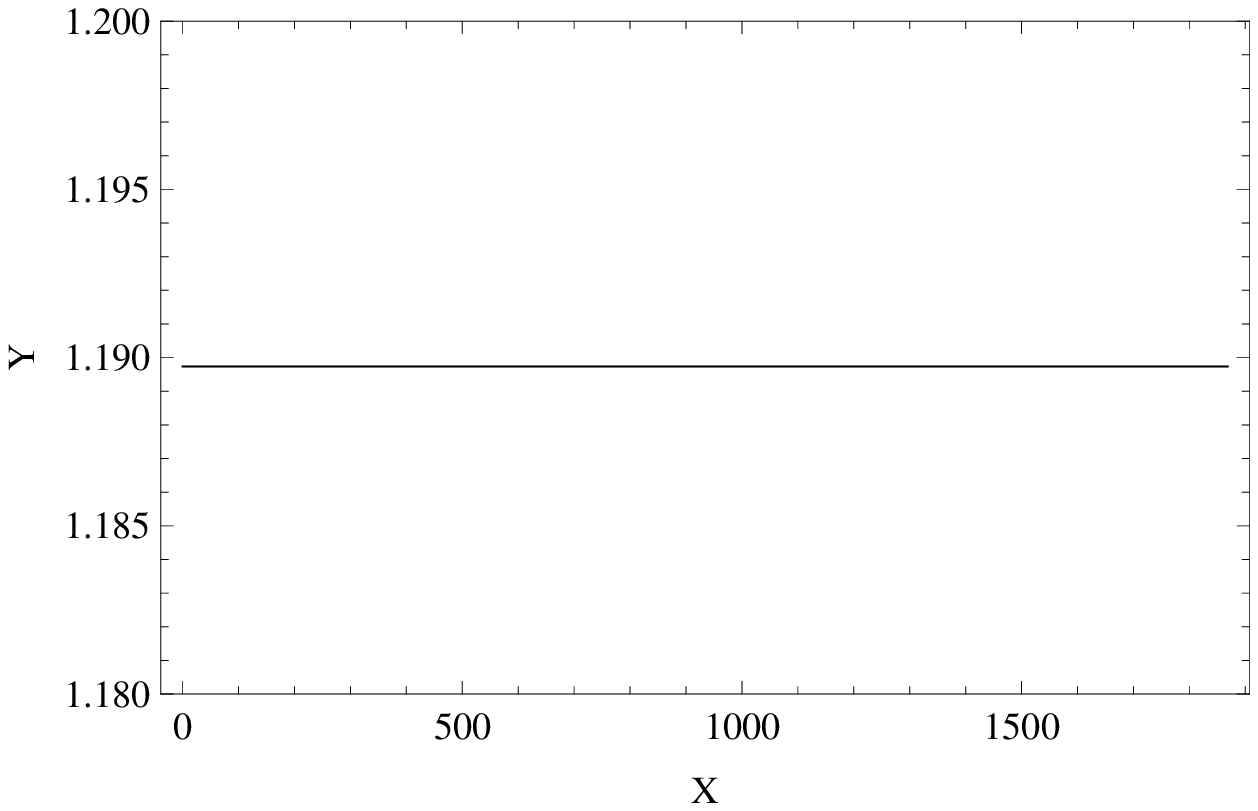}
\hglue-0.01cm
\includegraphics[width=4.8cm]{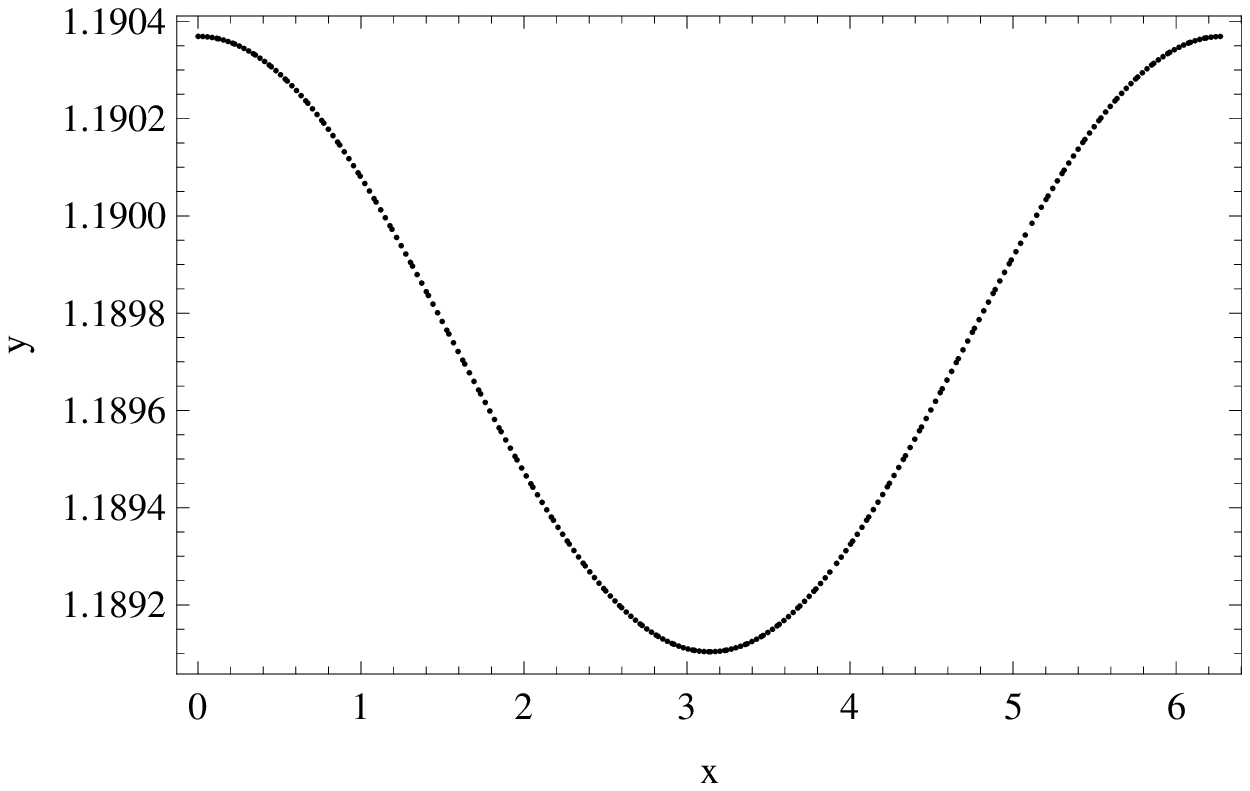}
\hglue-0.01cm
\includegraphics[width=5.1cm]{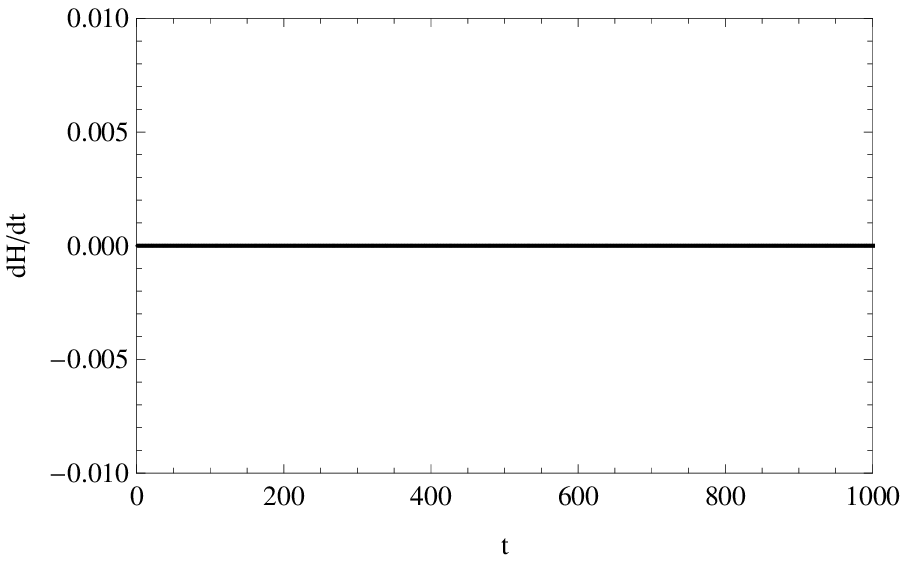}
\caption{Case $p_X=0$, $s\not=0$ associated to \equ{A1}
for $\varepsilon=10^{-3}$, $\mu=10^{-3}$ and for the initial conditions
$X(0)=0$, $Y(0)=1+6\sqrt{\varepsilon}$. Left: the lift of the normal form variables
$(X,Y)$ to the universal coverage. Middle: the trajectory in the original variables $(x,y)$.
Right: the variation of the derivative of the normalized Hamiltonian.}\label{case53}
\end{figure}

We consider an example for which the normal form equations provide $p_X=0$, but $s\neq0$.
To this end, we modify the conservative part, so that the actions do not contain a resonant term at first order:
\beqa{A1}
\dot x&=&y-\mu (\sin (x-t)+\sin (x)) \nonumber \\
\dot y&=&-\varepsilon (\sin (x-6t)+\sin (x)) -\mu  (y-\eta) \ .
\eeqa
The conservative transformation to second order is given by
\beqano
\psi_{10}(\tilde y,\tilde x,t)&=&\frac{\sin (\tilde x-6 t)}{\tilde y-6}+\frac{\sin (\tilde x)}{\tilde y} \nonumber \\
\psi_{20}(\tilde y,\tilde x,t)&=&-\frac{\sin (2 \tilde x)}{8 \tilde y^3}-\frac{\sin (2 \tilde x-6 t)}
{4 \tilde y^3-36 \tilde y^2+72 \tilde y}-\frac{\sin (2 \tilde x-12 t)}{8 (\tilde y-6)^3}
+\frac{\sin (6 t)}{72 \tilde y-12 \tilde y^2} \ .
\eeqano
The dissipative transformation to second order takes the form
\beqano
\beta_{01}(Y,X,t)&=&0\nonumber \\
\alpha_{01}(Y,X,t)&=&-\frac{\cos (X)}{Y} \nonumber \\
\beta_{11}(Y,X,t)&=&-\frac{\sin (2 X-7 t)}{4 Y^2-38 Y+84}-\frac{\sin (2 X-6 t)}{4 Y^2-36
   Y+72}-\frac{\sin (2 X)}{4 Y^2} \nonumber \\
&+&\frac{\sin (2 X-t)}{2 Y-4 Y^2}+\frac{\sin (t)}{2 Y}+\frac{\sin (5 t)}{10
   (Y-6)}+\frac{\sin (6 t)}{12 (Y-6)} \nonumber \\
\alpha_{11}(Y,X,t)&=&\frac{\cos (2 X)}{8 Y^3}+\frac{\cos (2 X-7 t)}{2 (7-2 Y)^2 (Y-6)}+\frac{\cos (2 X-6
   t)}{8 (Y-6) (Y-3)^2}\nonumber\\
&+&\frac{\cos (2 X-t)}{2 (1-2 Y)^2 Y}-\frac{(Y+2) \cos (t)}{2 Y^2}-\frac{(Y-16) \cos (5
   t)}{50 (Y-6)^2}-\frac{(Y-18) \cos (6 t)}{72 (Y-6)^2} \nonumber \\
\beta _{02}(Y,X,t)&=&0 \nonumber \\
\alpha _{02}(Y,X,t)&=&\frac{\sin (2 X)}{4 Y^2}+\frac{\sin (t)}{Y} \ .
\eeqano
The resulting normal form up to second order becomes
\beqano
\dot X&=&Y-\frac{\varepsilon ^2}{2 (Y-6)^3}-\frac{\varepsilon ^2}{2 Y^3}-\frac{\mu ^2}{2 Y} -\mu\sin (X-t) + O_3(\varepsilon,\mu) \nonumber \\
\dot Y &=& O_3(\varepsilon,\mu) \ ,
\eeqano
whereas the drift function is given by
$$
\eta(Y)=Y+{\varepsilon\over {2Y}}+O_3(\varepsilon,\mu)\ .
$$
The Hamiltonian function in normalized variables corresponding to $\mu=0$ in the extended phase space
turns out to be
\beqno
{\cal H}(Y,X,U,t)=\frac{Y^2}{2}+T+\frac{\left(Y^2-6 Y+18\right) \varepsilon ^2}{2 (Y-6)^2 Y^2} + O_3(\varepsilon,\mu)\ ;
\eeqno
the time derivative of the Hamiltonian under the dissipative flow becomes
\beqno
\frac {d {{\cal H}(Y,X,U,t)}} {d t} = O_3(\varepsilon,\mu) \ ,
\eeqno
which shows the preservation of the energy up to the third order.
Figure~\ref{case53} displays the behavior of the lift of $(X,Y)$ to the universal covering, the plot in the
original variables and the graph of the derivative of the Hamiltonian versus time. The result shows that the dynamics
takes place on an adiabatic quasi--periodic solution, which is consistent with the theoretical expectation.

\subsection{Exponential stability: case $p_X\not=0$, $s=0$}\label{EPXS2}
\begin{figure}[htp]
\includegraphics[width=4.8cm]{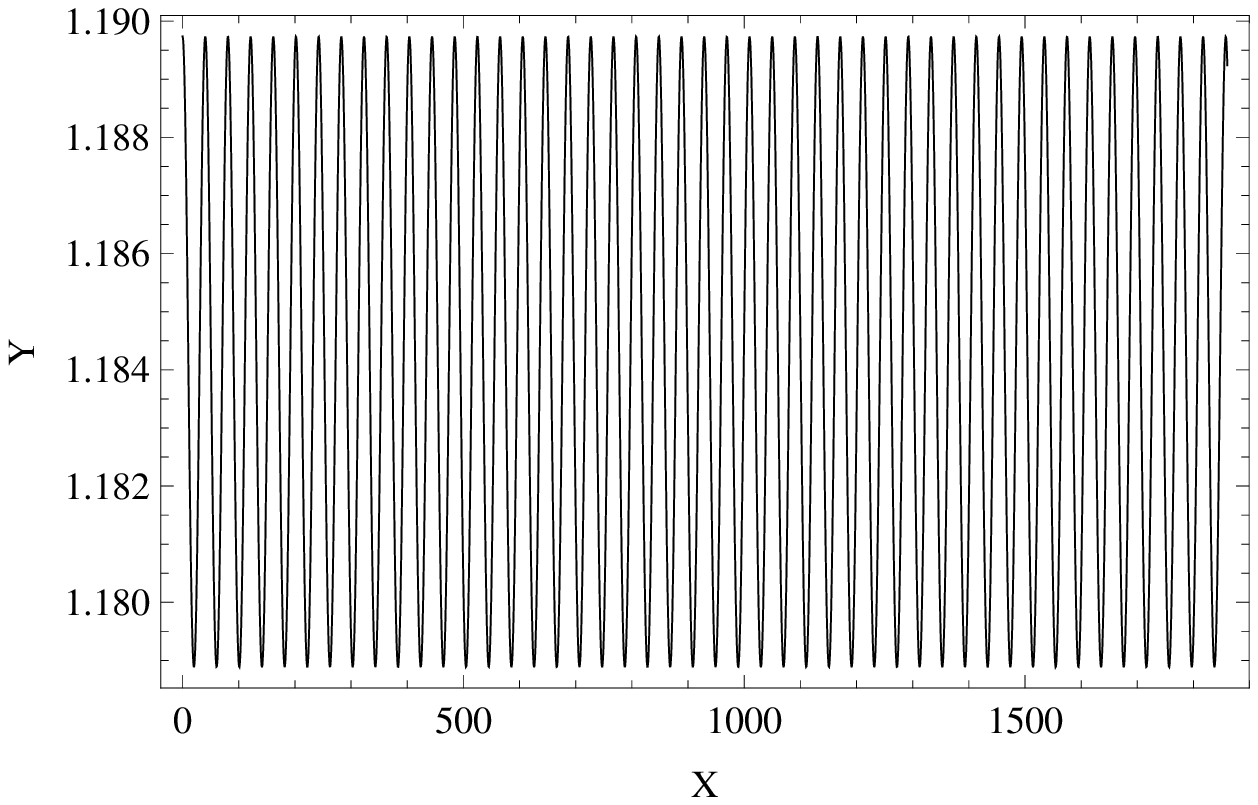}
\hglue-0.01cm
\includegraphics[width=4.8cm]{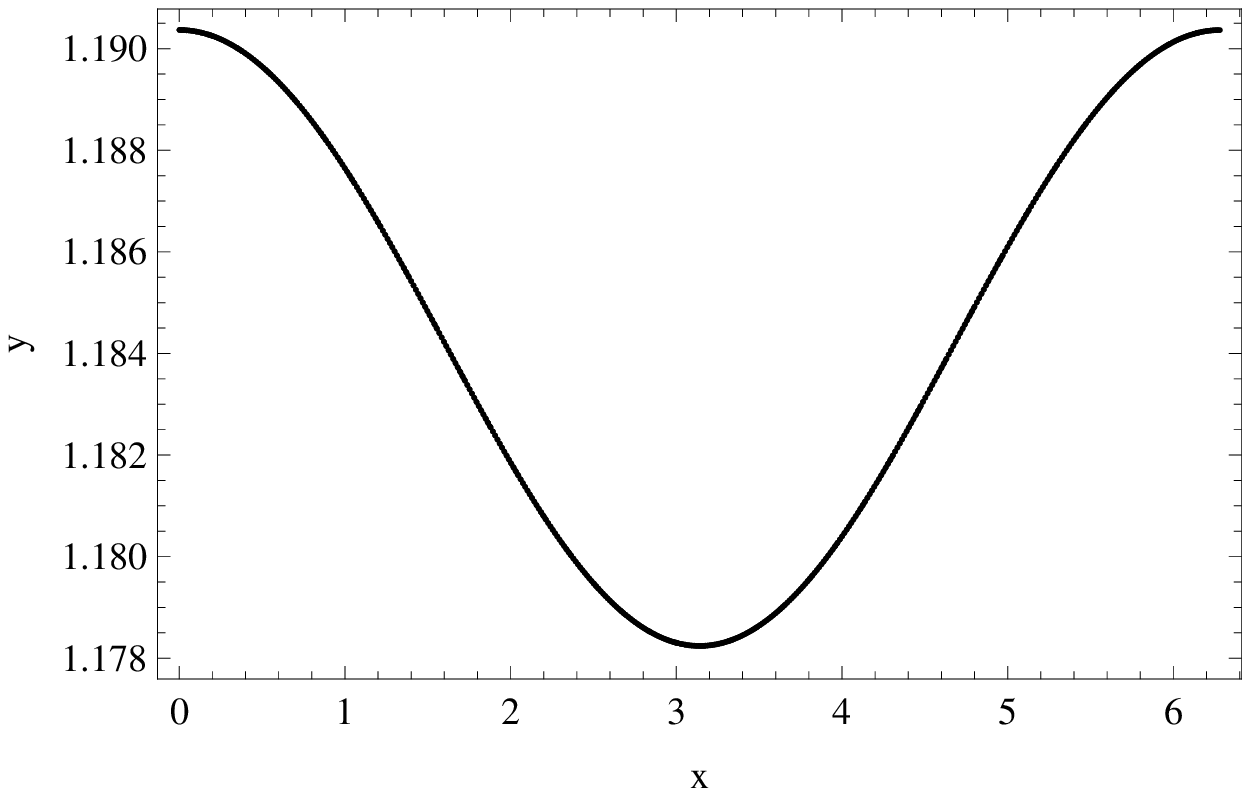}
\hglue-0.01cm
\includegraphics[width=5.1cm]{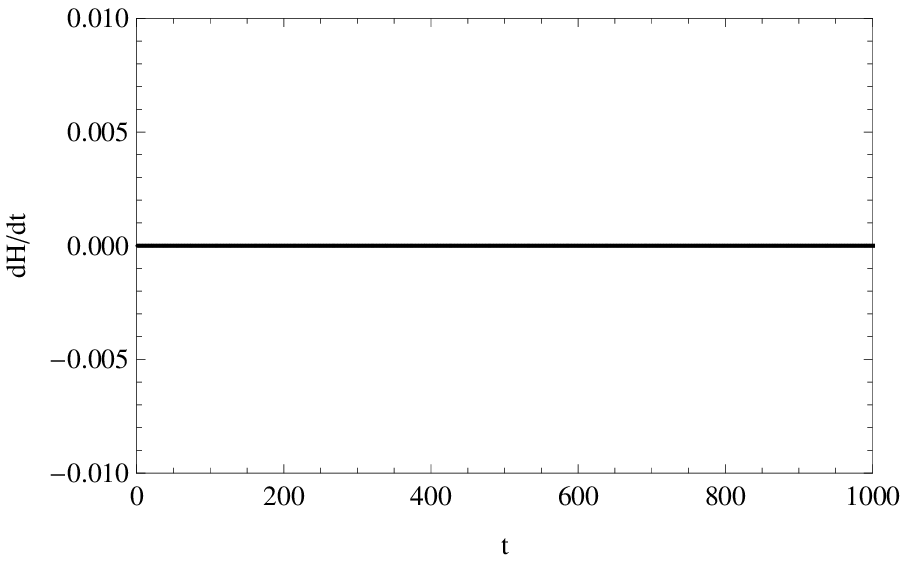}
\caption{Case $p_X\not=0$, $s=0$ associated to \equ{A2}
for $\varepsilon=10^{-3}$, $\mu=10^{-3}$ and for the initial conditions
$X(0)=0$, $Y(0)=1+6\sqrt{\varepsilon}$. Left: the lift of the normal form variables
$(X,Y)$ to the universal coverage. Middle: the trajectory in the original variables $(x,y)$.
Right: the variation of the derivative of the normalized Hamiltonian.}\label{case54}
\end{figure}

As an example which generates a normal form with $p_X\not=0$, $s=0$, we consider the differential equations
\beqa{A2}
\dot x&=&y-\mu\sin (6t) \nonumber \\
\dot y&=&-\varepsilon (\sin (x-t)+\sin (x)) -\mu  (y-\eta) \ .
\eeqa
The conservative transformation is given by
\beqano
\psi_{10}(\tilde y,\tilde x,t)&=&\frac{\sin (\tilde x)}{\tilde y} \nonumber \\
\psi_{20}(\tilde y,\tilde x,t)&=&-\frac{\sin (2 \tilde x)}{8 \tilde y^3}-\frac{\sin (2 \tilde x-t)}{2 \tilde y^2-4 \tilde y^3}-\frac{\sin (t)}{2 \tilde y^2} \ ,
\eeqano
while the dissipative transformation takes the form
\beqano
\beta_{01}(Y,X,t)&=&0\nonumber\\
\alpha_{01}(Y,X,t)&=&-{1\over 6}\cos(6t)\nonumber\\
\beta_{11}(Y,X,t)&=&-\frac{\sin (X-6 t)}{12 Y-2 Y^2}-\frac{\sin (X-7 t)}{84-12 Y}+\frac{\sin (X+5 t)}{12
   Y+60}-\frac{\sin (X+6 t)}{2 Y (Y+6)} \nonumber \\
\alpha_{11}(Y,X,t)&=&\frac{(3-Y) \cos (X-6 t)}{(Y-6)^2 Y^2}+\frac{(Y+3) \cos (X+6 t)}{Y^2
   (Y+6)^2}-\frac{\cos (X-7 t)}{12 (Y-7)^2}-\frac{\cos (X+5 t)}{12 (Y+5)^2}\nonumber \\
\beta_{02}(Y,X,t)&=&0 \nonumber \\
\alpha_{02}(Y,X,t)&=&0\ .
\eeqano
The normal form equations are given by
\beqano
\dot X&=&Y-\frac{\varepsilon ^2}{2 Y^3} + O_3(\varepsilon,\mu) \nonumber \\
\dot Y&=&-\varepsilon  \sin (X-t) + O_3(\varepsilon,\mu) \ ,
\eeqano
with the drift function provided by $\eta(Y,X,t)=Y-\frac{1}{144} \varepsilon  \mu  \sin (X-t)+O_3(\varepsilon,\mu)$.
Note that we produce linear conservative resonant terms in the actions,
but no resonant dissipative terms in the angles. The Hamiltonian function associated to the normal form equations
in the extended phase space becomes
\beqno
{\cal H}(Y,X,U,t)=\frac{Y^2}{2}+T+\frac{\varepsilon ^2}{4 Y^2}-\varepsilon  \cos (X-t) + O_3(\varepsilon,\mu)\ ,
\eeqno
while the time derivative of the Hamiltonian flow becomes;
\beqno
\frac{d {{\cal H}(Y,X,U,t)}}{d t} = O_3(\varepsilon,\mu) \ ,
\eeqno
yielding the preservation of the Hamiltonian up to the normalization order.
Figure~\ref{case54} shows the behavior of the lift of $(X,Y)$ to the universal covering, the graph in the
original variables and the plot of the derivative of the Hamiltonian versus time. Also in this case, the result shows that the dynamics
takes place on an adiabatic quasi--periodic solution, which is consistent with the theoretical expectation.

\section{Application of the stability estimates}\label{sec:APP}

In this Section we implement the Theorem to obtain estimates on the variation of the actions as given in Section~\ref{sec:stability}. Let us fix
the initial data as well as $r_0$, $s_0$ (and related domain's parameters), $K$, $\delta$.
We assume that the frequency satisfies \equ{res}, \equ{omega0} with $a$ determined by \equ{omega0}.
The smallness conditions on the parameters $\varepsilon$, $\mu$, say $\varepsilon\leq \varepsilon_0$, $\mu\leq\mu_0$, come from \equ{C1}, \equ{C2}, \equ{C2bis}, \equ{C4}, \equ{C5}, \equ{33ter}, \equ{C6},
\equ{cnew1}, \equ{C7}, \equ{C8}.

\noindent
We define the constants $\tilde C_G$ and $C_G$ as in \equ{cgtilde}, \equ{cg} and we set $C_Y\equiv C_G+\tilde C_G$.
We recall that $\tau_0$, $N$, $K$ are related by the expression
$$
\tau_0\equiv {N\over K}\ |\log \lambda|\ .
$$
From \equ{Pi} we determine $C_p$, while the constants $C_1$, $C_2$, $C_3$, $C_4$ are computed as in \equ{CC}.
Table~\ref{APPTAB} provides the main quantities involved in the Theorem through the application of a third order
normal form in the extended phase space. In particular, it provides the variation $\|\Delta Y\|$ of the normalized variables,
the variation $\|\Delta y\|\equiv\|y(t)-y(0)\|$ of the original variables and the stability time $T$, which perfectly
agrees with the theoretical result (linear or exponential stability time) of the Theorem.

\vskip.1in

\noindent
The results have been validated by a numerical integration of the equations of motion.
Due to computer limitations, for the cases described in Sections~\ref{EPXS1} and \ref{EPXS2}
we had to stop to a time at most equal to $10^8$. Up to such integration times the numerical
results are in full agreement with the analytical results.

\vskip.1in

\begin{table}
\centering
\(
\begin{array}{lllll}
\hline\hline
- & Sec.~\ref{EFLD} & Sec.~\ref{ESLD} & Sec.~\ref{EPXS1} & Sec.~\ref{EPXS2}  \\
\varepsilon_0 & 6.\cdot10^{-5} & 6.\cdot10^{-5} & 6.\cdot10^{-5} & 6.\cdot10^{-5} \\
\mu_0 & 6.\cdot10^{-5} & 6.\cdot10^{-5} & 6.\cdot10^{-5} & 1.9\cdot10^{-4} \\
\tau _0 & 1.458 & 1.458 & 1.458 & 1.285 \\
C_Y & 3.16\cdot10^{1} & 3.428\cdot10^{1} & 1.714\cdot10^{1} & 1.087 \\
C_p & 1.052  & 1.052 & 1.265 & 3.323\cdot10^{-1} \\
C_1 &    2.117\cdot10^{-3} & 2.233\cdot10^{-3} & 1.359\cdot10^{-3} & 2.158\cdot10^{-4} \\
C_2 &    5.056\cdot10^{-3} & 3.059\cdot10^{-5} & 0 & 0 \\
C_3 &    2.01 & 2.01     & 3.208\cdot10^{-5} & 2.01 \\
C_4 &    3.292\cdot10^{-5} & 3.292\cdot10^{-5} & 1.006\cdot10^{-5}. & 3.283\cdot10^{-6} \\
\|\Delta Y\| & 2.408\cdot10^{-2} & 2.408\cdot10^{-2} & 9.369\cdot10^{-3} & 2.408\cdot10^{-2} \\
\|\Delta y\| & 2.421\cdot10^{-2} & 2.421\cdot10^{-2} & 9.521\cdot10^{-3} & 2.421\cdot10^{-2} \\
T & 2.692\cdot10^{5} & 4.43\cdot10^{7} & 1.699\cdot10^{10} & 3.309\cdot10^{9} \\
\hline\hline
\end{array}
\)
\caption{The main quantities of the Theorem for the examples of Section~\ref{sec:examples} from
the remainder of a third order normal form.
The parameters and initial conditions for all columns are:
$x_0=0$, $y_0=1.01$, $r_0=0.05$, $\tilde r_0=4.9\cdot10^{-2}$, $\tilde r_0'=2.45\cdot10^{-2}$,
$R_0=2.4\cdot10^{-2}$, $s_0=0.1$, $\tilde s_0=5\cdot10^{-3}$, $S_0=2.5\cdot10^{-3}$,
$K=20$, $\delta=0.01$.}\label{APPTAB}
\end{table}


\newpage

\section{Appendix A}
We briefly review the conditions which must be satisfied by the parameters $\varepsilon$, $\mu$,
so that the transformation from original to intermediate variables, as well as that from intermediate
to final variables can be inverted; moreover, we provide conditions on the parameters so that the
non--resonance conditions in the intermediate and final variables are satisfied. Compare also with
\cite{cellho2010a} and \cite{GAL}.

\subsection{Inversion of the conservative transformation}
With reference to \equ{deltac2}, we invert the first transformation as
\beq{Ainv}
x=\tilde x+\Gamma^{(x,N)}(\tilde y,\tilde x,t)
\eeq
provided that
$$
70\, \|\psi_y^{(N)}\|_{\tilde r_0,s_0} e^{2s_0}\delta_0^{-1}<1\ ,
$$
with
$$
\|\Gamma^{(x,N)}\|_{\tilde r_0,\tilde s_0}\leq \|\psi_y^{(N)}\|_{\tilde r_0,s_0}\ .
$$
for $\tilde r_0<r_0$, $\delta_0<s_0$, $\tilde s_0\equiv s_0-\delta_0$.

\subsection{Non--resonance condition after the conservative transformation}
Taking into account \equ{omega0}, we want that the non--resonance condition is satisfied in the
intermediate variables, say for $a>0$:
\beq{omapp}
|\omega(\tilde y)\cdot k+m|>{a\over 2}\ ,\qquad |k|+|m|\leq K .
\eeq
The second of \equ{deltac2} can be inverted as
\beq{2b}
\tilde y=y+\varepsilon R^{(N)}(y,x,t)\ ,
\eeq
for a suitable function $R^{(N)}$ provided
$$
70\, \|\psi_x^{(N)}\|_{\tilde r_0,s_0}{1\over {\tilde r_0-\tilde r_0'}}<1\ ,
$$
for $\tilde r_0'<\tilde r_0$ with
$$
\varepsilon \|R^{(N)}\|_{\tilde r_0',s_0}\leq \|\psi_x^{(N)}\|_{\tilde r_0,s_0}\ .
$$
Then we have
\beqno
|\omega(\tilde y)\cdot k+m|
\geq {a\over 2}\ ,
\eeqno
if
$$
\varepsilon \leq {a\over {2K\|R^{(N)}\|_{\tilde r_0',s_0}\|\omega_y\|_{r_0}}}\ .
$$

\subsection{Inversion of the dissipative transformation}
With reference to \equ{B}, the first equation can be inverted provided
$$
70\,\|\alpha^{(N)}\|_{\tilde r_0,\tilde s_0}\ e^{2\tilde s_0}\tilde\delta_0^{-1}<1\ ,
$$
where $\tilde\delta_0<\tilde s_0$. Inverting the equation as
$$
\tilde x=X+ A^{(x,N)}(\tilde y,X,t)\ ,
$$
we have
$$
\|A^{(x,N)}\|_{\tilde r_0,\tilde s_0-\tilde\delta_0}\leq \|\alpha^{(N)}\|_{\tilde r_0,\tilde s_0}\ .
$$
Thus we invert the second of \equ{B} as
$$
\tilde y=Y+\Delta^{(y,N)}(Y,X,t)\ ,
$$
provided
$$
70\,\|A^{(y,N)}\|_{\tilde r_0,S_0}{1\over {\tilde r_0-R_0}}<1\ ,
$$
with $S_0<\tilde s_0-\tilde\delta_0$, $R_0<\tilde r_0$, being
$$
\|\Delta^{(y,N)}\|_{R_0,S_0}\leq \|A^{(y,N)}\|_{\tilde r_0,S_0}\ .
$$
Notice that $A^{(y,N)}$ can be bounded as
$$
\|A^{(y,N)}\|_{\tilde r_0,S_0}\leq\|\beta^{(N)}\|_{\tilde r_0,\tilde s_0}+
\|\beta^{(N)}_x\|_{\tilde r_0,\tilde s_0}\|A^{(x,N)}\|_{\tilde r_0,S_0}\ .
$$
Similar for the third equation in \equ{B}.

\subsection{Non--resonance condition after the dissipative transformation}
We now turn to the fulfillment of the non--resonant condition in the new set of variables
$$
|\omega(Y)\cdot k+m|>0\ ,\qquad |k|+|m|\leq K\ .
$$
Through the transformation
$$
Y=\tilde y+\beta^{(N)}(\tilde y,\tilde x,t;\varepsilon,\mu)
$$
and using \equ{omapp} one finds
\beqno
|\omega(Y)\cdot k+m|\geq{a\over 4}\ ,
\eeqno
provided that
$$
K\,\|\omega_y\|_{r_0} \|\beta^{(N)}\|_{\tilde r_0,\tilde s_0}<{a\over 4}\ .
$$

\vskip.2in

\section{Appendix B}
From properties of analytic functions one can prove the following result (see also \cite{cellho2010a}) on
the decay of the tail of the Fourier series.

\noindent
\bf Lemma B.1. \sl Let $f=f(y,x,t)$ be an analytic function on the domain
$C_{r_0}(A)\times C_{s_0}({\T}^{\ell+1})$. Let $f^{>K}(y,x,t)\equiv \sum_{(j,m)\in{\Z}^{\ell+1},
|j|+|m|>K} \hat f_{jm}(y)\ e^{i(j\cdot x+mt)}$ and let $0<\sigma_0<s_0$.
Then, there exists a constant $C_a\equiv C_a(\sigma_0,K)$, such that
\beq{uv}
\|f^{>K}\|_{r_0,s_0}\leq C_a \|f\|_{r_0,s_0+\sigma_0} e^{-(K+1)\sigma_0}\ ,
\eeq
with
\beq{Ca1}
C_a\equiv e^{(K+1){\sigma_0\over 2}}\left(\frac{1+e^{-{\sigma_0\over 2}}}{1-e^{-{\sigma_0\over 2}}}\right)^{\ell+1}\ .
\eeq
\rm

\vskip.3in

\noindent
\bf Acknowledgments. \rm
We are grateful to Luca Biasco, Enrico Valdinoci and Jean--Christophe Yoccoz for interesting discussions
and suggestions. We acknowledge the grants ASI
``Studi di Esplorazione del Sistema Solare" and PRIN 2007B3RBEY
``Dynamical Systems and Applications" of MIUR.

\end{document}